\documentclass[a4paper,10pt]{article}

\usepackage{paralist}
\usepackage{amsmath}
\usepackage[authoryear]{natbib}
\usepackage{tikz} 
\usetikzlibrary{calc,intersections}
\usepackage[top=3cm,bottom=3cm,left=4cm,right=4cm]{geometry}

\usepackage{hyperref}
\usepackage{subfigure}
\usepackage{MnSymbol}
\usepackage{booktabs}

\newcommand{\dd}{\; \mathrm{d}}

\newcommand{\dt}{\Delta t}

\newcommand{\paf}[2]{\frac{\partial #1}{\partial #2}}
\newcommand{\grad}{{\nabla}}
\newcommand{\free}{{\rm free}}
\newcommand{\jam}{{\rm jam}}

\newcommand{\rem}{{\rm remeshed}}
\newcommand{\init}{{\rm init}}

\newcommand{\stat}{{\rm stat}}
\newcommand{\dyn}{{\rm dyn}}
\newcommand{\old}{{\rm old}}
\newcommand{\new}{{\rm new}}
\newcommand{\mytriangle}{{\smalltriangleup}} 
\newcommand{\myvertex}{{\circ}}

\title{The Lagrangian coordinate system and what it means for two-dimensional crowd flow models}
\author{Femke van Wageningen-Kessels, Ludovic Leclercq, \\ Winnie Daamen and Serge P. Hoogendoorn}
\date{\today}

\begin{document}

\maketitle

\begin{abstract}
   A continuum crowd flow model is reformulated in the Lagrangian coordinate system. The system has proven to give computational advantages over the traditional Eulerian coordinate system for (one-dimensional) road traffic flow. Our extension of the model and simulation method to (two-dimensional) crowd flow paves the way to explore whether the advantages also hold in two dimensions. This paper provides a first exploration and shows that a model and simulation method for two-dimensional crowd flow can be developed in the Lagrangian coordinate system and that is leads to promising results.

\bigskip

\noindent
{\bf keywords} Crowd flow, Model, Simulation, Lagrangian coordinates
\end{abstract}


\section{Objective}
Crowd flow models are used to describe, understand and predict collective behaviour of crowds. Roughly, two types of models exist: microscopic models in which the movements of individual pedestrians are described and traced and macroscopic (or continuum) models which are the focus of this study. Continuum models describe the dynamics of crowds as a continuum flow, in terms of average speed, velocity and density \citep{Hughes2002_ContinuumPedestrian}. It has been shown that, in contrast to previous claims, continuum models can reproduce self-organisation and certain important dynamic phenomena such as lane formation and diagonal striping \citep{Hoogendoorn2014_MicroMacro_PhysicaA}. Furthermore, simulations based on continuum models have the potential to significantly reduce computation time whilst keeping high accuracy. This makes them useful for a larger range of applications, including real time state estimation and prediction for crowd management and optimisation of control strategies. For this, the solutions to the model equations need to be calculated using both fast and accurate computational methods.

The Lagrangian coordinate system has been applied to traffic flow models \citep{Leclercq2007}. \Citet{VanWageningen2010TRR,Yuan2012_IEEE} show that the system has many advantages over the traditional Eulerian coordinate system, including more accurate reproduction of shock waves and more efficient state estimation based on trajectory data. 

Our main contribution is the proposal of an extension of Lagrangian simulation methods previously applied for one-dimensional traffic flow models to two-dimensional crowd flow models. We introduce the reformulation of the two-dimensional continuum crowd flow model into Lagrangian coordinates (Section \ref{sec:LagrangianContinuumModel}), develop a numerical simulation method (Section \ref{sec:SimulationApproach}) and show that it leads to meaningful simulation results (Section \ref{sec:Simulation}). Open questions and future research directions are discussed in Section \ref{sec:Conclusion}.


\section{Lagrangian formulation of crowd flow model}
\label{sec:LagrangianContinuumModel}
The main idea of the Lagrangian (or moving) coordinate system is that the location $\vec x=\begin{pmatrix} x \\ y \end{pmatrix}$ of the $n$-th particle (i.e. a vehicle or a pedestrian) is calculated at time $t$: 
\begin{align}
   \vec x = \vec x (n,t)
\end{align}
This is in contrast to the Eulerian (fixed) coordinate system where the density at location $\vec x$ and time $t$ is calculated: 
\begin{align}
   \rho=\rho(\vec x,t)
\end{align}

In Eulerian coordinates, the two-dimensional conservation equation is well-known:
\begin{align}
   \paf{\rho}{t} + \paf{q_x}{x} + \paf{q_y}{y} = 0
\end{align}
with $\rho$ the density, and $q_x$ and $q_y$ the flow in $x$- and $y$-direction, respectively.
To derive this equation, Green's theorem can be applied. 

\subsection{Model in Lagrangian coordinates}
To derive the Lagrangian equivalent of the one-dimensional conservation equation in Eulerian coordinates, Green's theorem can be applied. However, for the two-dimensional case, this theorem can not be used to derive the conservation equation in the Lagrangian coordinate system. The main reason for this is that pedestrians can not be numbered consistently as they walk in a two-dimensional space. For example, they can overtake each other in an unpredictable way.


\begin{figure}
   \subfigure[Movement of one pedestrian. \label{fig:DerivationLag2D_1Ped}]{
   \begin{tikzpicture}[xscale=1.3]
   \coordinate (Cold) at (0,0);
   \coordinate (Cnew) at (2.5,1);
   \draw plot [smooth] coordinates {(Cold) (0.5,1) (1.2,1.5) (Cnew)};
   \fill (Cold) circle (2pt) node [right] {$\vec x_j(t_0)$};
   \fill (Cnew) circle (2pt) node [below] {$\vec x_j(t_1)$};
   \end{tikzpicture}
   }
   \hfill
   \subfigure[Movement of region with pedestrians. \label{fig:DerivationLag2D_Area}]{
   \begin{tikzpicture}[xscale=1.5]

   \coordinate (C00old) at (0,0);
   \coordinate (C05old) at (0.2,0.5);
   \coordinate (C10old) at (0.5,1);
   \coordinate (C15old) at (0.7,1.2);
   \coordinate (C20old) at (1.2,1.5);
   \coordinate (C25old) at (1.5,1.3);
   \coordinate (C30old) at (1.8,1.1);
   \coordinate (C35old) at (2.1,1.1);
   \coordinate (C40old) at (2.5,1);
   \coordinate (C45old) at (2.1,0.4);
   \coordinate (C50old) at (1.5,-0.2);
   \coordinate (C55old) at (0.6,-0.1);

   \coordinate (C00new) at ($(C00old)+(0,0.6)$);
   \coordinate (C05new) at ($(C05old)+(-0.2,0.8)$);
   \coordinate (C10new) at ($(C10old)+(-0.6,1)$);
   \coordinate (C15new) at ($(C15old)+(0,0.9)$);
   \coordinate (C20new) at ($(C20old)+(0.5,1)$);
   \coordinate (C25new) at ($(C25old)+(0.8,1)$);
   \coordinate (C30new) at ($(C30old)+(1.2,0.8)$);
   \coordinate (C35new) at ($(C35old)+(1.1,0.5)$);
   \coordinate (C40new) at ($(C40old)+(1,0.2)$);
   \coordinate (C45new) at ($(C45old)+(1.2,0.2)$);
   \coordinate (C50new) at ($(C50old)+(1.2,0.3)$);
   \coordinate (C55new) at ($(C55old)+(0.4,0.8)$);

   \begin{scope}[thick]
   \draw [fill=black!15,opacity=0.5] plot [smooth cycle] coordinates {(C00new) (C05new) (C10new) (C15new) (C20new) (C25new) (C30new) (C35new) (C40new) (C45new) (C50new) (C55new)};
   \draw[->] (C00old) -- (C00new); 
   \draw[->] (C10old) -- (C10new); 
   \draw[->] (C20old) -- (C20new); 
   \draw[->] (C30old) -- (C30new) node [left=5mm] {$A(t_1)$}; 
   \draw[->] (C40old) -- (C40new); 
   \draw[->] (C50old) node [above left, blue] {$A(t_0)$} -- (C50new); 
   \draw[->] (C55old) -- (C55new);
   \end{scope}

   \begin{scope}[blue]
   \draw plot [smooth cycle] coordinates {(C00old) (C05old) (C10old) (C15old) (C20old) (C25old) (C30old) (C35old) (C40old) (C45old) (C50old) (C55old)};
   \end{scope}

   \end{tikzpicture}
   }
\caption{Illustrations for the derivation of the discretised version of the two-dimensional conservation equation in Lagrangian coordinates.}
\label{fig:DerivationLag2D}
\end{figure}
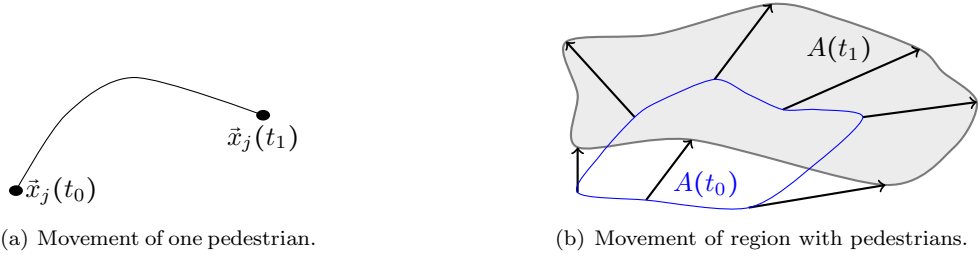

However, we can still derive a discretised version of the two-dimensional conservation equation in Lagrangian coordinates. It is based on its one-dimensional counterpart for road traffic \citep{Leclercq2007}. 
Therefore, we first consider the movement of one pedestrian $j$, which is at location $\vec x_j(t_0)$ at time $t=t_0$ and moves with velocity $\vec v_j(\vec x,t)$, see Figure \ref{fig:DerivationLag2D_1Ped}. At time $t=t_1>t_0$ its location will be:
\begin{align}
   \vec x_j(t_1) = \vec x_j(t_0) + \int_{\tau=t_0}^{t_1} v_j(\tau) \dd \tau
\label{eq:Movement1Ped}   
\end{align}
Secondly, we apply the continuum hypothesis (i.e. we do not consider individual pedestrians, but assume they are represented by a density $\rho(\vec x, t)$ at each location $\vec x$ and time $t$). Now we consider an region $i$ with $N_i(t) = \int_{\vec x \in \Omega_i(t)} \rho(\vec x, t) \dd \vec x$ pedestrians in it. The area of this region equals:
\begin{align}
   A_i(t) = \int_{\vec x \in \Omega_i(t)} 1 \dd \vec x = \oint_{\vec x \in \Gamma_i(t)} \vec x \dd \vec x
\label{eq:ComputeArea}   
\end{align}
with $\Gamma$ the boundary of the region.
In the Lagrangian coordinate system, the region moves with the pedestrians, and thus the number of pedestrians does not change ($N(t_1)=N(t_0)$). However, the region and its boundary move if speed is nonzero ($\Gamma_i(t_1) \neq \Gamma_i(t_0)$).
Using both \eqref{eq:Movement1Ped} and \eqref{eq:ComputeArea}, we can compute the area $A(t_1)$ at time $t=t_1>t_0$, see also Figure \ref{fig:DerivationLag2D_Area}:
\begin{align}
   A_i(t_1) &= \oint_{\vec x \in \Gamma_i(t_1)} \vec x \dd \vec x
              = \oint_{\vec x \in \Gamma_i(t_0)} \vec x + \int_{\tau=t}^{t_1} v(\vec x,\tau) \dd \tau \dd \vec x
              \nonumber \\
              &= A_i(t_0) + \oint_{\vec x \in \Gamma_i(t_0)} \int_{\tau=t_0}^{t_1} v(\vec x,\tau) \dd \tau \dd \vec x
\label{eq:ComputeAreaNew}   
\end{align}
The above concept of change in area occupied by a group of pedestrians according to the velocity of the pedestrians at the boundary, forms the basis of the numerical computation method, as outlined in Section \ref{sec:SimulationApproach}.

\subsection{Velocity: speed and route choice}
\label{sec:VelocityModel}
To complete the model, we also need to describe the velocity $\vec v(\vec x, t)$.
We propose a simple velocity model, which may be applied either in Lagrangian or in Eulerian coordinates. The velocity is calculated using a simplified version of the approach described by \citet{Hoogendoorn2014_MicroMacro_PhysicaA,Hoogendoorn2015_ISTTT_LocalRouteChoice}:
\begin{align}
   \vec v(\vec x,t) = V(\rho(\vec x, t)) 
      \frac{ \vec e_\stat (\vec x, t) + \beta_\dyn \vec e_\dyn (\rho(\vec x, t)) }{ \| \vec e_\stat (\vec x, t) + \beta_\dyn \vec e_\dyn (\rho(\vec x, t)) \| }
\label{eq:RouteChoice}      
\end{align}
with $V(\rho)$ the fundamental relation, $\vec e_\stat$ the static route choice component, $\vec e_\dyn$ the dynamic (local, state dependent) route choice component and $\beta_\dyn$ its weight (a model parameter). 
We chose simple models for speed and route choice to illustrate the numerical method, without too much focus on the continuous model. Other, more realistic models for speed and route choice can be implemented as well.

\subsubsection{Speed}
The speed only depends on density, following a linear fundamental relation:
\begin{align}
   V(\rho) = v_\free \left( 1 - \frac{\rho}{\rho_\jam} \right)
\label{eq:FD}
\end{align}
with $v_\free$ the free flow speed and $\rho_\jam$ the jam density.

\subsubsection{Route choice: static and dynamic}
The static route choice component $\vec e_\stat$ depends on the layout of the facility and possibly other factors. It can, for example reflect the shortest route to the destination. In the simulations (Section \ref{sec:Simulation}) we apply three different static route choices, which are described in more detail there.

For the dynamic route choice we apply a simple model for pedestrians avoiding high density regions:
\begin{align}
   \vec e_\dyn = - \grad \rho
\label{eq:DynRouteChoice}   
\end{align}
The dynamic route choice component steers pedestrians in the direction with the steepest negative density gradient.

\section{Simulation approach}
\label{sec:SimulationApproach}
The main idea of the simulation approach using the Lagrangian coordinate system is that a moving grid is applied. In particular, the grid moves with the same velocity as the pedestrians. In the following we discuss how the grid is constructed, how the velocity of its vertices is computed and how the grid is moved each time step $\dt$. Furthermore, we discuss when and how the grid is reconstructed.

\subsection{The grid}
\begin{figure}
\centering
\begin{tikzpicture}[xscale=0.4,yscale=0.4*0.75*sqrt(3)]

   \clip (4.5,2.5) rectangle (17.5,11.5);


   
   \coordinate (V02_00) at (2,0);
   \coordinate (V06_00) at (6,0);
   \coordinate (V10_00) at (10,0);
   \coordinate (V14_00) at (14,0);
   \coordinate (V18_00) at (18,0);

   \coordinate (V00_03) at (0,3);
   \coordinate (V04_03) at (4,3);
   \coordinate (V08_03) at (8,3);
   \coordinate (V12_03) at (12,3);
   \coordinate (V16_03) at (16,3);
   \coordinate (V20_03) at (20,3);

   \coordinate (V02_06) at (2,6);
   \coordinate (V06_06) at (6,6);
   \coordinate (V10_06) at (10,6);
   \coordinate (V14_06) at (14,6);
   \coordinate (V18_06) at (18,6);

   \coordinate (V00_09) at (0,9);
   \coordinate (V04_09) at (4,9);
   \coordinate (V08_09) at (8,9);
   \coordinate (V12_09) at (12,9);
   \coordinate (V16_09) at (16,9);
   \coordinate (V20_09) at (20,9);

   \coordinate (V02_12) at (2,12);
   \coordinate (V06_12) at (6,12);
   \coordinate (V10_12) at (10,12);
   \coordinate (V14_12) at (14,12);
   \coordinate (V18_12) at (18,12);

   \coordinate (C00_01) at (0,1);
   \coordinate (C04_01) at (4,1);
   \coordinate (C08_01) at (8,1);
   \coordinate (C12_01) at (12,1);
   \coordinate (C16_01) at (16,1);
   \coordinate (C20_01) at (20,1);
   
   \coordinate (C02_02) at (2,2);
   \coordinate (C06_02) at (6,2);
   \coordinate (C10_02) at (10,2);
   \coordinate (C14_02) at (14,2);
   \coordinate (C18_02) at (18,2);
   
   \coordinate (C00_05) at (0,5);
   \coordinate (C04_05) at (4,5);
   \coordinate (C08_05) at (8,5);
   \coordinate (C12_05) at (12,5);
   \coordinate (C16_05) at (16,5);
   \coordinate (C20_05) at (20,5);
   
   \coordinate (C02_04) at (2,4);
   \coordinate (C06_04) at (6,4);
   \coordinate (C10_04) at (10,4);
   \coordinate (C14_04) at (14,4);
   \coordinate (C18_04) at (18,4);
   
   \coordinate (C00_07) at (0,7);
   \coordinate (C04_07) at (4,7);
   \coordinate (C08_07) at (8,7);
   \coordinate (C12_07) at (12,7);
   \coordinate (C16_07) at (16,7);
   \coordinate (C20_07) at (20,7);
   
   \coordinate (C02_08) at (2,8);
   \coordinate (C06_08) at (6,8);
   \coordinate (C10_08) at (10,8);
   \coordinate (C14_08) at (14,8);
   \coordinate (C18_08) at (18,8);
   
   \coordinate (C00_11) at (0,11);
   \coordinate (C04_11) at (4,11);
   \coordinate (C08_11) at (8,11);
   \coordinate (C12_11) at (12,11);
   \coordinate (C16_11) at (16,11);
   \coordinate (C20_11) at (20,11);
   
   \coordinate (C02_10) at (2,10);
   \coordinate (C06_10) at (6,10);
   \coordinate (C10_10) at (10,10);
   \coordinate (C14_10) at (14,10);
   \coordinate (C18_10) at (18,10);

   \coordinate (Voffset) at (0,-0.04);


   \draw (V02_00) -- (V06_00) -- (V10_00) -- (V14_00) -- (V18_00);
   \draw (V00_03) -- (V02_00) -- (V04_03) -- (V06_00) -- (V08_03) -- (V10_00) -- (V12_03) -- (V14_00) -- (V16_03) -- (V18_00) -- (V20_03);
   \draw (V00_03) -- (V04_03) -- (V08_03) -- (V12_03) -- (V16_03) -- (V20_03);
   \draw (V00_03) -- (V02_06) -- (V04_03) -- (V06_06) -- (V08_03) -- (V10_06) -- (V12_03) -- (V14_06) -- (V16_03) -- (V18_06) -- (V20_03);
   \draw (V02_06) -- (V06_06) -- (V10_06) -- (V14_06) -- (V18_06);
   \draw (V00_09) -- (V02_06) -- (V04_09) -- (V06_06) -- (V08_09) -- (V10_06) -- (V12_09) -- (V14_06) -- (V16_09) -- (V18_06) -- (V20_09);
   \draw (V00_09) -- (V04_09) -- (V08_09) -- (V12_09) -- (V16_09) -- (V20_09);
   \draw (V00_09) -- (V02_12) -- (V04_09) -- (V06_12) -- (V08_09) -- (V10_12) -- (V12_09) -- (V14_12) -- (V16_09) -- (V18_12) -- (V20_09);
   \draw (V02_12) -- (V06_12) -- (V10_12) -- (V14_12) -- (V18_12);

   \begin{scope}[font={\large}]
   \node at ($(V02_00)+(Voffset)$) {$\myvertex$};
   \node at ($(V06_00)+(Voffset)$) {$\myvertex$};
   \node at ($(V10_00)+(Voffset)$) {$\myvertex$};
   \node at ($(V14_00)+(Voffset)$) {$\myvertex$};
   \node at ($(V18_00)+(Voffset)$) {$\myvertex$};

   \node at ($(V00_03)+(Voffset)$) {$\myvertex$};
   \node at ($(V04_03)+(Voffset)$) {$\myvertex$};
   \node at ($(V08_03)+(Voffset)$) {$\myvertex$};
   \node at ($(V12_03)+(Voffset)$) {$\myvertex$};
   \node at ($(V16_03)+(Voffset)$) {$\myvertex$};
   \node at ($(V20_03)+(Voffset)$) {$\myvertex$};

   \node at ($(V02_06)+(Voffset)$) {$\myvertex$};
   \node at ($(V06_06)+(Voffset)$) {$\myvertex$};
   \node at ($(V10_06)+(Voffset)$) {$\myvertex$};
   \node at ($(V14_06)+(Voffset)$) {$\myvertex$};
   \node at ($(V18_06)+(Voffset)$) {$\myvertex$};

   \node at ($(V00_09)+(Voffset)$) {$\myvertex$};
   \node at ($(V04_09)+(Voffset)$) {$\myvertex$};
   \node at ($(V08_09)+(Voffset)$) {$\myvertex$};
   \node at ($(V12_09)+(Voffset)$) {$\myvertex$};
   \node at ($(V16_09)+(Voffset)$) {$\myvertex$};
   \node at ($(V20_09)+(Voffset)$) {$\myvertex$};

   \node at ($(V02_12)+(Voffset)$) {$\myvertex$};
   \node at ($(V06_12)+(Voffset)$) {$\myvertex$};
   \node at ($(V10_12)+(Voffset)$) {$\myvertex$};
   \node at ($(V14_12)+(Voffset)$) {$\myvertex$};
   \node at ($(V18_12)+(Voffset)$) {$\myvertex$};
   
   \node at ($(C00_01)+(Voffset)$) {$\mytriangle$};
   \node at ($(C04_01)+(Voffset)$) {$\mytriangle$};
   \node at ($(C08_01)+(Voffset)$) {$\mytriangle$};
   \node at ($(C12_01)+(Voffset)$) {$\mytriangle$};
   \node at ($(C16_01)+(Voffset)$) {$\mytriangle$};
   \node at ($(C20_01)+(Voffset)$) {$\mytriangle$};

   \node at ($(C02_02)+(Voffset)$) {$\mytriangle$};
   \node at ($(C06_02)+(Voffset)$) {$\mytriangle$};
   \node at ($(C10_02)+(Voffset)$) {$\mytriangle$};
   \node at ($(C14_02)+(Voffset)$) {$\mytriangle$};
   \node at ($(C18_02)+(Voffset)$) {$\mytriangle$};

   \node at ($(C00_05)+(Voffset)$) {$\mytriangle$};
   \node at ($(C04_05)+(Voffset)$) {$\mytriangle$};
   \node at ($(C08_05)+(Voffset)$) {$\mytriangle$};
   \node at ($(C12_05)+(Voffset)$) {$\mytriangle$};
   \node at ($(C16_05)+(Voffset)$) {$\mytriangle$};
   \node at ($(C20_05)+(Voffset)$) {$\mytriangle$};

   \node at ($(C02_04)+(Voffset)$) {$\mytriangle$};
   \node at ($(C06_04)+(Voffset)$) {$\mytriangle$};
   \node at ($(C10_04)+(Voffset)$) {$\mytriangle$};
   \node at ($(C14_04)+(Voffset)$) {$\mytriangle$};
   \node at ($(C18_04)+(Voffset)$) {$\mytriangle$};

   \node at ($(C00_07)+(Voffset)$) {$\mytriangle$};
   \node at ($(C04_07)+(Voffset)$) {$\mytriangle$};
   \node at ($(C08_07)+(Voffset)$) {$\mytriangle$};
   \node at ($(C12_07)+(Voffset)$) {$\mytriangle$};
   \node at ($(C16_07)+(Voffset)$) {$\mytriangle$};
   \node at ($(C20_07)+(Voffset)$) {$\mytriangle$};

   \node at ($(C02_08)+(Voffset)$) {$\mytriangle$};
   \node at ($(C06_08)+(Voffset)$) {$\mytriangle$};
   \node at ($(C10_08)+(Voffset)$) {$\mytriangle$};
   \node at ($(C14_08)+(Voffset)$) {$\mytriangle$};
   \node at ($(C18_08)+(Voffset)$) {$\mytriangle$};

   \node at ($(C00_11)+(Voffset)$) {$\mytriangle$};
   \node at ($(C04_11)+(Voffset)$) {$\mytriangle$};
   \node at ($(C08_11)+(Voffset)$) {$\mytriangle$};
   \node at ($(C12_11)+(Voffset)$) {$\mytriangle$};
   \node at ($(C16_11)+(Voffset)$) {$\mytriangle$};
   \node at ($(C20_11)+(Voffset)$) {$\mytriangle$};

   \node at ($(C02_10)+(Voffset)$) {$\mytriangle$};
   \node at ($(C06_10)+(Voffset)$) {$\mytriangle$};
   \node at ($(C10_10)+(Voffset)$) {$\mytriangle$};
   \node at ($(C14_10)+(Voffset)$) {$\mytriangle$};
   \node at ($(C18_10)+(Voffset)$) {$\mytriangle$};

   \end{scope}

\end{tikzpicture}
\caption{A staggered triangular grid. Density $\rho_\mytriangle$ is calculated at the cell centres, velocity $\vec v_\myvertex$ is calculated at the vertices.}
\label{fig:StaggeredGrid}
\end{figure}
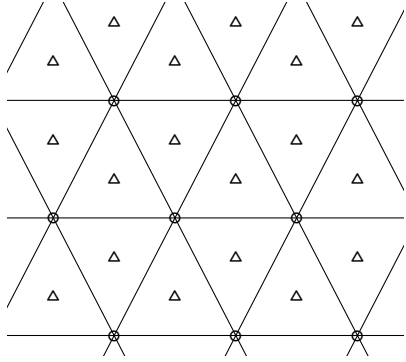

We assume triangular grid cells, but the methods can be extended to other polygon-shaped grid cells. The grid is staggered, meaning that some variables are calculated at the cell centres, others at the vertices, see Figure \ref{fig:StaggeredGrid}. We do not consider any variables calculated at the edges. The cell centres are defined as the centre of mass of the cells. Variables calculated at the cell centres are indicated with a subscripted triangle ($\mytriangle$). Variables calculated at the vertex are indicated with a subscripted circle ($\myvertex$). The main variable to be calculated at the cells centres is density ($\rho_\mytriangle$). The main variable to be calculated at the vertices is velocity ($\vec v_\myvertex$).

\subsection{Velocity of grid cell vertices}
We know the location of the grid cells and their vertices and the number of pedestrians within each cell from the previous time step. From this we approximate the velocity of the grid vertices, as outlined below.

\subsubsection{Density in grid cells}
Each time step, the average density in the grid cell is calculated:
\begin{align}
   \rho_\mytriangle = \frac{N_\mytriangle}{A_\mytriangle}
\end{align}
and serves as an approximation of the density at the grid cell centre.
$N_\mytriangle$ is the number of pedestrians in the grid cell (remains constant over time). $A_\mytriangle$ is the area of the grid cell (may change each time step), see also Figure \ref{fig:NumberingVertex}. In general it is computed as in \eqref{eq:ComputeArea}, but because the grid cell is triangular we can use:
\begin{align}
   A_\triangle = {\textstyle \frac12} \left|
      ( x_{j_1} - x_{j_3} )( y_{j_2} - y_{j_1} ) -  ( x_{j_1} - x_{j_2} )( y_{j_3} - y_{j_1} )
   \right|
\end{align} 
with the locations of the vertices $\{j_1,j_2,j_3\} \in J$ of the triangle:
\begin{align}
   \vec x_{\myvertex,j} = \begin{pmatrix} x_j \\ y_j \end{pmatrix}
\end{align}

\subsubsection{Speed at vertices}
To obtain the speeds from the fundamental relation, the densities at the vertices are approximated as a weighted average of the densities of the surrounding cells:
\begin{align}
   \rho_\myvertex = \sum_{i\in I} \alpha_i \rho_{\mytriangle,i} 
\end{align}
with $i \in I$ the indices of all surrounding cells.
The weights $\alpha_i$ are proportional to the proximity of the cell centre to the vertex (see also Figure \ref{fig:NumberingCell}):
\begin{align}
   \alpha_i = \frac{1/ \| \vec x_{\mytriangle,i} - \vec x_\myvertex \|}{ \sum_{i^*\in I} \left( 1 / \| \vec x_{\mytriangle,i^*} - \vec x_\myvertex \| \right)}
\end{align}
with $\vec x_\myvertex$ the location of the vertex. The cell centres $\vec x_{\mytriangle,i}$ are calculated as the centre of mass of the cells:
\begin{align}
   \vec x_{\mytriangle,i} = {\textstyle \frac13} \sum_{j \in J_i} \vec x_{\myvertex,j}
\end{align}
with $j \in J_i$ the indices of the vertices of cell $i$.
Subsequently, the speed of the vertices are approximated using the fundamental relation $v_\myvertex=V(\rho_\myvertex)$, see equation \eqref{eq:FD}.

\subsubsection{Route choice at vertices}
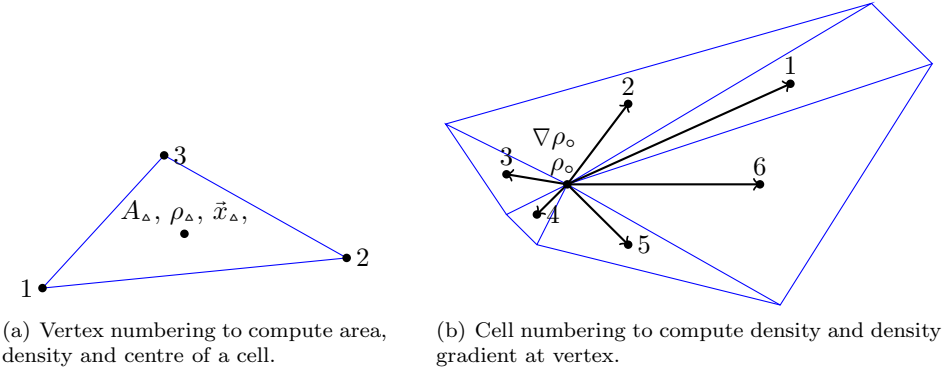
\begin{figure}[t]
\hfill
   \subfigure[Vertex numbering to compute area, density and centre of a cell.\label{fig:NumberingVertex}]{
		\begin{tikzpicture}[scale = 0.8]
		   \coordinate (V1) at (-2.5,0);
		   \coordinate (V2) at (2.5,0.5);   
		   \coordinate (V3) at (-0.5,2.2);
		   
		   \draw[blue] (V1) -- (V2) -- (V3) -- cycle;
		   
		   \fill (V1) circle (2pt) node [left]  {$1$};
		   \fill (V2) circle (2pt) node [right] {$2$};
		   \fill (V3) circle (2pt) node [right] {$3$};
		   
		   \coordinate (M2) at ($(V1)!+0.5!(V2)$);
		   \coordinate (CC) at ($(V3)!+2/3!(M2)$);
		   
		   \fill (CC) circle (2pt) node [above] {$A_{\mytriangle}$, $\rho_{\mytriangle}$, $\vec x_{\mytriangle}$, };
		\end{tikzpicture}
	}
\hfill
   \subfigure[Cell numbering to compute density and density gradient at vertex.\label{fig:NumberingCell}]{
\begin{tikzpicture}[scale = 0.8]

	\pgfdeclarelayer{background}
	\pgfdeclarelayer{foreground}
	\pgfsetlayers{background,main,foreground}

   
   \coordinate (P0) at (4,3);
   \coordinate (P1) at (10,5);
   \coordinate (P2) at (9,6);
   \coordinate (P3) at (2,4);
   \coordinate (P4) at (3,2.5);
   \coordinate (P5) at (3.5,2);
   \coordinate (P6) at (7.5,1);
   
   \begin{scope}[blue]
   \draw (P0) -- (P1) -- (P2);
   \draw (P0) -- (P2) -- (P3);
   \draw (P0) -- (P3) -- (P4);
   \draw (P0) -- (P4) -- (P5);
   \draw (P0) -- (P5) -- (P6);
   \draw (P0) -- (P6) -- (P1);
   \end{scope}
   
	 \coordinate (M1) at ($(P1)!+0.5!(P2)$);
	 \coordinate (M2) at ($(P2)!+0.5!(P3)$);
	 \coordinate (M3) at ($(P3)!+0.5!(P4)$);
	 \coordinate (M4) at ($(P4)!+0.5!(P5)$);
	 \coordinate (M5) at ($(P5)!+0.5!(P6)$);
	 \coordinate (M6) at ($(P6)!+0.5!(P1)$);

    \coordinate (C1) at ($(P0)!+2/3!(M1)$);
    \coordinate (C2) at ($(P0)!+2/3!(M2)$);
    \coordinate (C3) at ($(P0)!+2/3!(M3)$);
    \coordinate (C4) at ($(P0)!+2/3!(M4)$);
    \coordinate (C5) at ($(P0)!+2/3!(M5)$);
    \coordinate (C6) at ($(P0)!+2/3!(M6)$);
   
    \fill (P0) circle (2pt) node [above,text width=10mm,text centered] {$\grad \rho_\myvertex \quad$ \vspace*{-5mm} \\ $\rho_\myvertex \;$};
   
    \fill (C1) circle (2pt) node [above] {$1$};
    \fill (C2) circle (2pt) node [above] {$2$};
    \fill (C3) circle (2pt) node [above] {$3$};
    \fill (C4) circle (2pt) node [right] {$4$};
    \fill (C5) circle (2pt) node [right] {$5$};
    \fill (C6) circle (2pt) node [above] {$6$};

    \begin{scope}[->,thick]
    \draw (P0) -- (C1);
    \draw (P0) -- (C2);
    \draw (P0) -- (C3);
    \draw (P0) -- (C4);
    \draw (P0) -- (C5);
    \draw (P0) -- (C6);
    \end{scope}

\end{tikzpicture}
}
\hfill
\caption{Examples of vertex and cell numbering. This is used to compute the area, density and centre of a cell \subref{fig:NumberingVertex} and to compute the density and density gradient at a vertex \subref{fig:NumberingCell}. In the example in \subref{fig:NumberingCell}, density $\rho_1$ in cell $i=1$ will get a low weight $\alpha_1$ because its centre is far from the vertex for which we calculate density $\rho_\myvertex$ and density gradient $\grad \rho_\myvertex$. Density $\rho_4$ in cell $i=4$ will get a high weight $\alpha_4$ because its centre is close to the vertex.}
\label{fig:DensityAtVertex}
\end{figure}

For the dynamic route choice component we need the density gradient at the vertices. They are approximated as the weighted average of the gradients over each of the lines between the vertex and the centres of the surrounding cells (see also Figure \ref{fig:NumberingCell}):
\begin{align}
   \grad \rho_\myvertex 
   = \sum_{i\in I} \alpha_i \grad \rho_{\mytriangle,i}
   = \sum_{i\in I} \alpha_i 
     \left( \rho_{\mytriangle,i} - \rho_{\myvertex} \right)
     \begin{pmatrix}
	     \frac{1}{x_{\mytriangle,i} - x_\myvertex}
	     \smallskip \\
	     \frac{1}{y_{\mytriangle,i} - y_\myvertex}
     \end{pmatrix}
\end{align}
We use $\vec e_{\dyn,\myvertex} = - \grad \rho_\myvertex$ (see equation \eqref{eq:DynRouteChoice}) to calculate the dynamic route choice component.

The static route choice component $\vec e_\stat$ is supposed to be given at any location and time.

\subsubsection{Velocity of vertices}
We can now approximate the velocity of the vertex (see equation \eqref{eq:RouteChoice}):
\begin{align}
   \vec v_\myvertex 
   = v_\myvertex
      \frac{ \vec e_{\stat,\myvertex} - \beta_\dyn \grad \rho_\myvertex }{ \| \vec e_{\stat,\myvertex} - \beta_\dyn \grad \rho_\myvertex \| }
\label{eq:VelocityDiscrete}      
\end{align}

\subsection{Moving of grid cells}
\begin{figure}
\centering
\begin{tikzpicture}[xscale=0.4,yscale=0.4*0.75*sqrt(3)]

   \clip (4.5,2.5) rectangle (17.5,11.5);


   \coordinate (V02_00) at (2,0);
   \coordinate (V06_00) at (6,0);
   \coordinate (V10_00) at (10,0);
   \coordinate (V14_00) at (14,0);
   \coordinate (V18_00) at (18,0);

   \coordinate (V00_03) at (0,3);
   \coordinate (V04_03) at (4,3);
   \coordinate (V08_03) at (8,3);
   \coordinate (V12_03) at (12,3);
   \coordinate (V16_03) at (16,3);
   \coordinate (V20_03) at (20,3);

   \coordinate (V02_06) at (2,6);
   \coordinate (V06_06) at (6,6);
   \coordinate (V10_06) at (10,6);
   \coordinate (V14_06) at (14,6);
   \coordinate (V18_06) at (18,6);

   \coordinate (V00_09) at (0,9);
   \coordinate (V04_09) at (4,9);
   \coordinate (V08_09) at (8,9);
   \coordinate (V12_09) at (12,9);
   \coordinate (V16_09) at (16,9);
   \coordinate (V20_09) at (20,9);

   \coordinate (V02_12) at (2,12);
   \coordinate (V06_12) at (6,12);
   \coordinate (V10_12) at (10,12);
   \coordinate (V14_12) at (14,12);
   \coordinate (V18_12) at (18,12);

   \coordinate (W02_00) at ($(2,0)+(rand*1.25+0.75,rand*1.25+0.75)$);
   \coordinate (W06_00) at ($(6,0)+(rand*1.25+0.75,rand*1.25+0.75)$);
   \coordinate (W10_00) at ($(10,0)+(rand*1.25+0.75,rand*1.25+0.75)$);
   \coordinate (W14_00) at ($(14,0)+(rand*1.25+0.75,rand*1.25+0.75)$);
   \coordinate (W18_00) at ($(18,0)+(rand*1.25+0.75,rand*1.25+0.75)$);

   \coordinate (W00_03) at ($(0,3)+(rand*1.25+0.75,rand*1.25+0.75)$);
   \coordinate (W04_03) at ($(4,3)+(rand*1.25+0.75,rand*1.25+0.75)$);
   \coordinate (W08_03) at ($(8,3)+(1.5,1.1)$);
   \coordinate (W12_03) at ($(12,3)+(rand*1.25+0.75,rand*1.25+0.75)$);
   \coordinate (W16_03) at ($(16,3)+(rand*1.25+0.75,rand*1.25+0.75)$);
   \coordinate (W20_03) at ($(20,3)+(rand*1.25+0.75,rand*1.25+0.75)$);

   \coordinate (W02_06) at ($(2,6)+(rand*1.25+0.75,rand*1.25+0.75)$);
   \coordinate (W06_06) at ($(6,6)+(1,0.8)$);
   \coordinate (W10_06) at ($(10,6)+(0.9,1.4)$);
   \coordinate (W14_06) at ($(14,6)+(rand*1.25+0.75,rand*1.25+0.75)$);
   \coordinate (W18_06) at ($(18,6)+(rand*1.25+0.75,rand*1.25+0.75)$);

   \coordinate (W00_09) at ($(0,9)+(rand*1.25+0.75,rand*1.25+0.75)$);
   \coordinate (W04_09) at ($(4,9)+(rand*1.25+0.75,rand*1.25+0.75)$);
   \coordinate (W08_09) at ($(8,9)+(rand*1.25+0.75,rand*1.25+0.75)$);
   \coordinate (W12_09) at ($(12,9)+(rand*1.25+0.75,rand*1.25+0.75)$);
   \coordinate (W16_09) at ($(16,9)+(rand*1.25+0.75,rand*1.25+0.75)$);
   \coordinate (W20_09) at ($(20,9)+(rand*1.25+0.75,rand*1.25+0.75)$);

   \coordinate (W02_12) at ($(2,12)+(rand*1.25+0.75,rand*1.25+0.75)$);
   \coordinate (W06_12) at ($(6,12)+(rand*1.25+0.75,rand*1.25+0.75)$);
   \coordinate (W10_12) at ($(10,12)+(rand*1.25+0.75,rand*1.25+0.75)$);
   \coordinate (W14_12) at ($(14,12)+(rand*1.25+0.75,rand*1.25+0.75)$);
   \coordinate (W18_12) at ($(18,12)+(rand*1.25+0.75,rand*1.25+0.75)$);

   \coordinate (C00_01) at (0,1);
   \coordinate (C04_01) at (4,1);
   \coordinate (C08_01) at (8,1);
   \coordinate (C12_01) at (12,1);
   \coordinate (C16_01) at (16,1);
   \coordinate (C20_01) at (20,1);
   
   \coordinate (C02_02) at (2,2);
   \coordinate (C06_02) at (6,2);
   \coordinate (C10_02) at (10,2);
   \coordinate (C14_02) at (14,2);
   \coordinate (C18_02) at (18,2);
   
   \coordinate (C00_05) at (0,5);
   \coordinate (C04_05) at (4,5);
   \coordinate (C08_05) at (8,5);
   \coordinate (C12_05) at (12,5);
   \coordinate (C16_05) at (16,5);
   \coordinate (C20_05) at (20,5);
   
   \coordinate (C02_04) at (2,4);
   \coordinate (C06_04) at (6,4);
   \coordinate (C10_04) at (10,4);
   \coordinate (C14_04) at (14,4);
   \coordinate (C18_04) at (18,4);
   
   \coordinate (C00_07) at (0,7);
   \coordinate (C04_07) at (4,7);
   \coordinate (C08_07) at (8,7);
   \coordinate (C12_07) at (12,7);
   \coordinate (C16_07) at (16,7);
   \coordinate (C20_07) at (20,7);
   
   \coordinate (C02_08) at (2,8);
   \coordinate (C06_08) at (6,8);
   \coordinate (C10_08) at (10,8);
   \coordinate (C14_08) at (14,8);
   \coordinate (C18_08) at (18,8);
   
   \coordinate (C00_11) at (0,11);
   \coordinate (C04_11) at (4,11);
   \coordinate (C08_11) at (8,11);
   \coordinate (C12_11) at (12,11);
   \coordinate (C16_11) at (16,11);
   \coordinate (C20_11) at (20,11);
   
   \coordinate (C02_10) at (2,10);
   \coordinate (C06_10) at (6,10);
   \coordinate (C10_10) at (10,10);
   \coordinate (C14_10) at (14,10);
   \coordinate (C18_10) at (18,10);

   \coordinate (Voffset) at (0,-0.08);


   \fill [blue,opacity=0.3] (V08_03) -- (V10_06) -- (V06_06) -- cycle;
   \fill [black,opacity=0.3] (W08_03) -- (W10_06) -- (W06_06) -- cycle;

   \begin{scope}[blue,very thin]
   \draw (V02_00) -- (V06_00) -- (V10_00) -- (V14_00) -- (V18_00);
   \draw (V00_03) -- (V02_00) -- (V04_03) -- (V06_00) -- (V08_03) -- (V10_00) -- (V12_03) -- (V14_00) -- (V16_03) -- (V18_00) -- (V20_03);
   \draw (V00_03) -- (V04_03) -- (V08_03) -- (V12_03) -- (V16_03) -- (V20_03);
   \draw (V00_03) -- (V02_06) -- (V04_03) -- (V06_06) -- (V08_03) -- (V10_06) -- (V12_03) -- (V14_06) -- (V16_03) -- (V18_06) -- (V20_03);
   \draw (V02_06) -- (V06_06) -- (V10_06) -- (V14_06) -- (V18_06);
   \draw (V00_09) -- (V02_06) -- (V04_09) -- (V06_06) -- (V08_09) -- (V10_06) -- (V12_09) -- (V14_06) -- (V16_09) -- (V18_06) -- (V20_09);
   \draw (V00_09) -- (V04_09) -- (V08_09) -- (V12_09) -- (V16_09) -- (V20_09);
   \draw (V00_09) -- (V02_12) -- (V04_09) -- (V06_12) -- (V08_09) -- (V10_12) -- (V12_09) -- (V14_12) -- (V16_09) -- (V18_12) -- (V20_09);
   \draw (V02_12) -- (V06_12) -- (V10_12) -- (V14_12) -- (V18_12);
   \end{scope}

   \draw (W02_00) -- (W06_00) -- (W10_00) -- (W14_00) -- (W18_00);
   \draw (W00_03) -- (W02_00) -- (W04_03) -- (W06_00) -- (W08_03) -- (W10_00) -- (W12_03) -- (W14_00) -- (W16_03) -- (W18_00) -- (W20_03);
   \draw (W00_03) -- (W04_03) -- (W08_03) -- (W12_03) -- (W16_03) -- (W20_03);
   \draw (W00_03) -- (W02_06) -- (W04_03) -- (W06_06) -- (W08_03) -- (W10_06) -- (W12_03) -- (W14_06) -- (W16_03) -- (W18_06) -- (W20_03);
   \draw (W02_06) -- (W06_06) -- (W10_06) -- (W14_06) -- (W18_06);
   \draw (W00_09) -- (W02_06) -- (W04_09) -- (W06_06) -- (W08_09) -- (W10_06) -- (W12_09) -- (W14_06) -- (W16_09) -- (W18_06) -- (W20_09);
   \draw (W00_09) -- (W04_09) -- (W08_09) -- (W12_09) -- (W16_09) -- (W20_09);
   \draw (W00_09) -- (W02_12) -- (W04_09) -- (W06_12) -- (W08_09) -- (W10_12) -- (W12_09) -- (W14_12) -- (W16_09) -- (W18_12) -- (W20_09);
   \draw (W02_12) -- (W06_12) -- (W10_12) -- (W14_12) -- (W18_12);

   \begin{scope}[red,->,thick]
   \draw (V08_03) -- (W08_03);
   \draw (V06_06) -- (W06_06);
   \draw (V10_06) -- (W10_06);
   \end{scope}

\end{tikzpicture}
\caption{Example of movement of the grid during one time step. The old grid is shown in thin blue lines. The new grid in thick black lines. One grid cell is highlighted and the figure shows how it moves and deforms, following the velocities of its corners (red arrows).}
\label{fig:MoveGrid}
\end{figure}
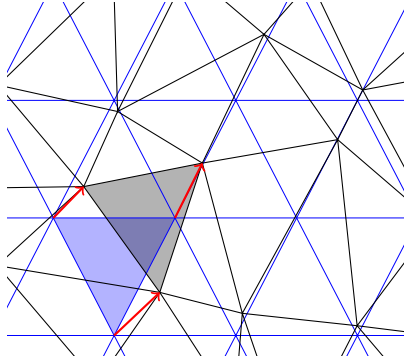

Once the velocity of the vertices is known, they can be moved. Therefore, we assume that the velocity remains as computed in \eqref{eq:VelocityDiscrete} during the time step. In every time step, the vertices are moved as follows:
\begin{align}
   \vec x_\myvertex^{\new} = \vec x_\myvertex^{\old} + \dt \vec v_\myvertex
\end{align}
This is illustrated in Figure \ref{fig:MoveGrid}.
We assume that the velocity of the edges is a linear interpolation of the velocities of the neighbouring vertices. Therefore, edges remain straight and the cell remains in the shape of a triangle.

\subsection{Remeshing of grid}
The velocities of the vertices of an individual grid cell, will usually have different velocities. Therefore, the grid cell will change its shape. If this deforms the grid cells `too much' it may be necessary to remesh. 

\subsubsection{Conditions for remeshing}
If either of the following holds, a new mesh is formed:
\begin{compactenum}
\item the density in at least one of the grid cells exceeds jam density,
\item in at least one of the grid cells the order of the vertices has changed: e.g. if vertices were first numbered in counter clockwise direction, they are now numbered in clockwise direction or vice versa,
\item the area of the triangle has become so small that 
\begin{align}
   \frac{A_{i,\new}}{A_{i,\init}} < \alpha
\end{align}
with $A_{i,\new}$ the new area of cell $i$, $A_{i,\init}$ its initial area and $0<\alpha<1$ a remeshing parameter.
\end{compactenum}
In the second case, a change of ordering indicates that a triangular grid cell has `flipped'. In the third case, a small fraction indicates that the grid cell has become too `stretched'. It may not be necessary for stability of the numerical method in this case, but we expect that it will improve accuracy.

\subsubsection{Remeshing procedure}
In the remeshing process a new mesh is formed and the pedestrians in the old grid cells, are assigned to new grid cells in such a way that the density at each location $\vec x$ remains (approximately) the same and the total number of pedestrians in the computational domain does not change (conservation of pedestrians).

\begin{figure}
\centering
\begin{tikzpicture}

	\pgfdeclarelayer{background}
	\pgfdeclarelayer{foreground}
	\pgfsetlayers{background,main,foreground}

   
   \coordinate (P0) at (4,3);
   \coordinate (P1) at (7,5);
   \coordinate (P2) at (6,6);
   \coordinate (P3) at (1.5,5);
   \coordinate (P4) at (0,2.5);
   \coordinate (P5) at (2.5,0.5);
   \coordinate (P6) at (6,1);
   \coordinate (P7) at (8.5,1.5);

   \coordinate (Q1) at (2,1.5);
   \coordinate (Q2) at (4.5,5.4528);
   \coordinate (Q3) at (7,1.5);
   
   \path[name path=l61]  (P6) -- (P1);
   
   \path[name path=l1]   (P0) -- (P1);
   \path[name path=l2]   (P0) -- (P2);
   \path[name path=l3]   (P0) -- (P3);
   \path[name path=l4]   (P0) -- (P4);
   \path[name path=l5]   (P0) -- (P5);
   \path[name path=l6]   (P0) -- (P6);
   
   \path[name path=lq12]   (Q1) -- (Q2);
   \path[name path=lq23]   (Q2) -- (Q3);
   \path[name path=lq31]   (Q3) -- (Q1);
 
   \path [name intersections={of = l1 and lq23}];
   \coordinate (p1)  at (intersection-1);
   \path [name intersections={of = l2 and lq23}];
   \coordinate (p2)  at (intersection-1);
   \path [name intersections={of = l3 and lq12}];
   \coordinate (p3)  at (intersection-1);
   \path [name intersections={of = l4 and lq12}];
   \coordinate (p4)  at (intersection-1);
   \path [name intersections={of = l5 and lq31}];
   \coordinate (p5)  at (intersection-1);
   \path [name intersections={of = l6 and lq31}];
   \coordinate (p6)  at (intersection-1);
   \path [name intersections={of = l61 and lq31}];
   \coordinate (p61a)  at (intersection-1);
   \path [name intersections={of = l61 and lq23}];
   \coordinate (p61b)  at (intersection-1);
   
   \begin{scope}[opacity=0.2,thick]
   \filldraw[draw=black,fill=yellow]  (P0) -- (P1) -- (P2);
   \filldraw[draw=black,fill=green]   (P0) -- (P2) -- (P3);
   \filldraw[draw=black,fill=cyan]    (P0) -- (P3) -- (P4);
   \filldraw[draw=black,fill=blue]    (P0) -- (P4) -- (P5);
   \filldraw[draw=black,fill=magenta] (P0) -- (P5) -- (P6);
   \filldraw[draw=black,fill=red]     (P0) -- (P6) -- (P1);
   \filldraw[draw=black,fill=orange]  (P6) -- (P7) -- (P1);
   \end{scope}
   
   \begin{scope}[opacity=0.5]
   \fill[yellow]  (P0) -- (p1) -- (p2);
   \fill[green]   (P0) -- (p2) -- (Q2) -- (p3);
   \fill[cyan]    (P0) -- (p3) -- (p4);
   \fill[blue]    (P0) -- (p4) -- (Q1) -- (p5);
   \fill[magenta] (P0) -- (p5) -- (p6);
   \fill[red]     (P0) -- (p6) -- (p61a) -- (p61b) -- (p1);
   \fill[orange]  (p61a) -- (Q3) -- (p61b);
   \end{scope}
   
   \draw[thick]   (Q1) -- (Q2) -- (Q3) -- cycle;
   
   \node at (6,5)      {$A_1$};   \node at (5,4)      {$A_1^*$};
   \node at (3.5,5)    {$A_2$};   \node at (4,4)    {$A_2^*$};
   \node at (2,3.5)    {$A_3$};   \node at (3.35,3.15)    {$A_3^*$};
   \node at (1.5,2)    {$A_4$};   \node at (3,2.25)    {$A_4^*$};
   \node at (4,1)      {$A_5$};   \node at (4.25,2)      {$A_5^*$};
   \node at (6.25,4)   {$A_6$};   \node at (5.25,2.75)   {$A_6^*$};
   \node at (7.25,2.5) {$A_7$};   \node at (6.5,1.75) {$A_7^*$};

\end{tikzpicture}
\caption{Example of remeshing procedure. The old mesh consists of (lightly coloured) triangles $A_i$. One triangle of the new mesh is drawn with thick black lines. It overlaps with 7 triangles of the old mesh. The intersections of the old and new triangles are brightly coloured and are indicated with $A_i^*$.}
\label{fig:Remesh}
\end{figure}
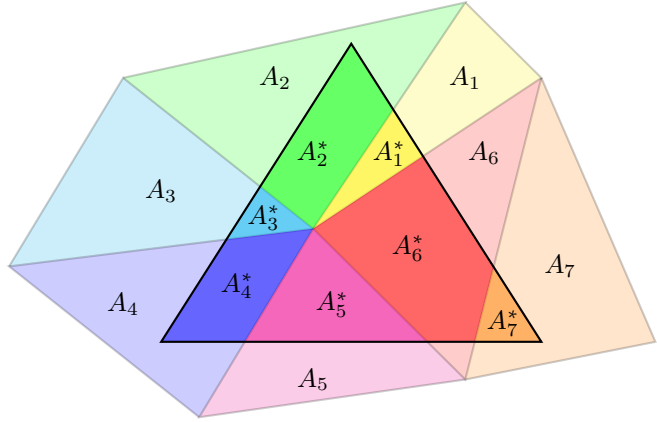

In our approach, the mesh after remeshing is the same as the initial mesh at the beginning of the simulation. This is a regular mesh consisting of equilateral triangles of the same size, just as in Figure \ref{fig:StaggeredGrid}.
Figure \ref{fig:Remesh} shows an example of the remeshing procedure. For the remeshing, density in the $j$-th triangle of the new mesh $\rho_j^\rem$ is calculated as follows:
\begin{align}
   \rho_j^\rem 
   = \frac{N_j^\rem}{A_j^\rem}
   = \frac{1}{A_j^\rem}\sum_i \frac{A_i^*}{A_i} N_i
\label{eq:Remesh}
\end{align}
with $A_i$ the area of the $i$-th triangle of the old mesh, $N_i$ the number of pedestrians in that triangle and $A_i^*$ the area of the intersection between the $i$-th triangle of the old mesh and the $j$-th triangle in the new mesh.

\section{Simulation and results}
\label{sec:Simulation}
We perform simulations with the newly developed method. We consider three test cases for which only the static route choice and the initial density are different. All other settings are equal for all test cases. We study the plausibility of the results and identify directions for improvements.

\subsection{Test set up}
We use a computational domain that is big enough to neglect the influence of the boundaries. 
The initial mesh consists of equilateral triangles. The parameters of the model and simulation are given in Table \ref{tab:Parameters}. We chose the test cases such that we have a variety of situations, with high and low densities, increasing and decreasing over time and with changes in walking direction of the pedestrians. 

\begin{table}[p]
\caption{Parameters of model and simulation}
\label{tab:Parameters}
{\bf Model parameters} 
\\
\begin{tabular}{l l r@{.}l l}
\hline 
free flow speed              & $v_\free$    &1&3 & m/s \\
jam density                  & $\rho_\jam$  &5&4 & m$^{-2}$ \\
weight dynamic route choice  & $\beta_\dyn$ &0&5 & \\
amplitude parameter (case 2) & $a$    &\multicolumn{2}{c}{$\pi / 2$} &   \\
frequency parameter (case 2) & $b$    &\multicolumn{2}{c}{$1 / (2\pi)$} & m$^{-1}$ \\
location of spiral centre (case 3)  & $(x_c,y_c)$ &\multicolumn{2}{c}{(60,60)}  \\
attraction of centre parameter (case 3) & $b$    &0&2 &  \\
\hline
\end{tabular}
\medskip \\
{\bf Initial condition parameters}
\\
\begin{tabular}{l l l l l}
\hline 
                            &                    &case 1 and 2& case 3     \\
location of centre of group & $(X_0,Y_0)$        &$(-30,0)$   &$(-30,0)$     & m \\
radius of group             & $R_0$              &$20$        &$10$          & m \\
\hline
\end{tabular}
\medskip \\
{\bf Numerical parameters}
\\
\begin{tabular}{l l l l l}
\hline 
area of initial triangles      & $A$    &56.9      & m$^2$ \\
initial $x$-range of domain    & X      &$[0,120]$ & m \\   
initial $y$-range of domain    & Y      &$[0,120]$ & m \\   
time step size                 & $\dt$  &1         & s \\
remeshing parameter         & $\alpha$  &0.01      &  \\
\hline
\end{tabular}
\end{table}

\subsubsection{Case 1 (`straight')}
In test case 1 (`straight') the static route choice directs the pedestrians from left to right, as shown in Figure \ref{fig:StaticStraight}. It is described by:
\begin{align}
   \vec e = 
   \begin{pmatrix}
	   1 \\ 0
   \end{pmatrix} 
\end{align}
As a result, streamlines are described by
\begin{align}
   y(x) = y_0
\end{align}
with $y_0 = y(0)$ the intersection of the streamline with the $y$-axis.

Initially, the pedestrians are located in a group centred at $(X_0,Y_0)$, see Figure \ref{fig:Initial}. The density is high at the centre of the group, and decreases linearly with distance to the centre, until it reaches zero at distance $R_0$:
\begin{align}
   \rho_0 = \rho(x,y,0) = 
	\begin{cases} \rho_\jam (1 - r/20) & \mbox{if $r<R_0$} \\
                 0                    & \mbox{if $r\geq R_0$} 
   \end{cases}
\label{eq:IC_Straight_ZigZag}
\end{align}
with $r=\sqrt{(x-X_0)^2+(y-Y_0)^2}$ the distance to the centre of the group. Parameter values are given in Table \ref{tab:Parameters}.

\subsubsection{Case 2 (`zig-zag')}
In test case 2 (`zig-zag') the static route choice directs the pedestrians on a zig-zag path, as shown in Figure \ref{fig:StaticZigZag}. It is described by:
\begin{align}
   \vec e = \frac{1}{ \sqrt{ a^2 + \sin^2 (bx) } } \begin{pmatrix}
   a \\ \sin \left( b x \right) \end{pmatrix}
\end{align}
with $a$ and $b$ route choice parameters, with values as in Table \ref{tab:Parameters}. As a result, streamlines are described by
\begin{align}
   y(x) = \frac{1 - \cos \left( bx \right)}{ab} + y_0
\end{align}
with $y_0 = y(0)$ the intersection of the streamline with the $y$-axis.

Initially, the pedestrians are located just as in test case 1, see \eqref{eq:IC_Straight_ZigZag} and Figure \ref{fig:Initial}.

\subsubsection{Case 3 (`spiral')}
In test case 3 (`spiral') the static route choice directs the pedestrians in a spiral-shape towards the centre of the spiral $(x_c,y_c)$, as shown in Figure \ref{fig:StaticSpiral}. It is described by:
\begin{align}
   \vec e = \frac{ 1 }{ \sqrt{ \left( b\tilde x-r\tilde y \right)^2 + \left( b\tilde y+r\tilde x \right)^2 }  }
   \begin{pmatrix}
	   b\tilde x-r\tilde y \\ b\tilde y+r\tilde x
   \end{pmatrix} 
\end{align}
with $\tilde x=x-x_c$ the distance to the centre of the spiral in $x$-direction, $\tilde y=y-y_c$ the distance to the centre of the spiral in $y$-direction and $r=\sqrt{\tilde x^2+\tilde y ^2}$ the absolute distance to the centre of the spiral. $b$ is a parameter indicating the strength of the attraction of the centre, with value as in Table \ref{tab:Parameters}.

As a result, streamlines are as follows, see also Figure \ref{fig:StaticSpiral}. Consider a pedestrian originating at $(x_0,y_0)$. Its original distance from the centre of the spiral is$r_0=\sqrt{(x_0-x_c)^2+(y_0-y_c)^2}$. The angle $\theta(x,y)$ is defined as the angle of the line through the centre $(x_c,y_c)$ and the position of the pedestrian $(x,y)$ with the line parallel to the $x$-axis and through the centre $(x_c,y_c)$. When the pedestrian walks, the angle $\theta$ increases and its distance to the centre $r$ decreases as follows:
\begin{align}
   r = \max(0, r_0 - b (\theta - \theta_0) )
\end{align}
with $\theta_0=\theta(x_0,y_0)$ the angle at the original position of the pedestrian. 

\begin{figure}[p]
\centering
   \begin{tikzpicture}[scale=1]
	   
	   \begin{scope}[thick]
	   \shade [inner color=gray,outer color=white,draw=black] (-1.75,0) circle (1.05);
	   \fill (-1.75,0) circle (2pt) node [above=2mm] {$\rho=\rho_0(r)$} node [above=10mm] {$\rho=0$};
	   \draw[->] (-1.75,0) -- +(-40:1.05); 
	   \node at (-0.7,0) [below left] {$R_0$};

	   \draw (-1.75,0) -- (-1.75,-1.45) node [below] {$X_0$};
	   \draw (-1.75,0) -- (-3.2,0) node [left] {$Y_0$};

		\end{scope}

	   \draw [->] (-3.5,-1.25) -- (1,-1.25) node [at end,below left] {$x$};
	   \draw [->] (-3,-1.75) -- (-3,1.75) node [at end,above left,sloped] {$y$};
	   
   \end{tikzpicture}
\caption{Initial condition.}
\label{fig:Initial}
\end{figure}
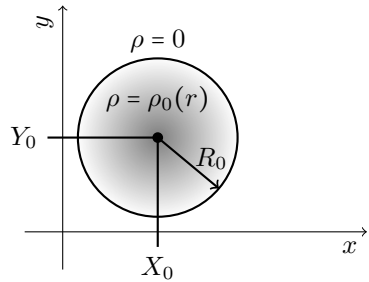

For the simulation, initially, the pedestrians are located in a group centred at $(X_0,Y_0)$, just as in the other 2 test cases, see Figure \ref{fig:Initial}. However, in this case, the density is constant over the region with radius $R_0$: 
\begin{align}
   \rho_0 = \rho(x,y,0) = 
	\begin{cases} \frac12 \rho_\jam    & \mbox{if $r<R_0$} \\
                 0                    & \mbox{if $r\geq R_0$} 
   \end{cases}
\label{eq:IC_Spiral}
\end{align}
with $r=\sqrt{(x-X_0)^2+(y-Y_0)^2}$ the distance to the centre of the group. Parameter values are given in Table \ref{tab:Parameters}.

\begin{figure}[p]
\subfigure[Straight (case 1).\label{fig:StaticStraight}]{
   \begin{tikzpicture}[scale=0.75]
	   	   \clip (-1,-1) rectangle ++(3.5,6);

	   \foreach \oneoverab in {0.75}{
		   \begin{scope}[gray]
		   \foreach \yinit in {-4,-3.5,...,6}{
		      \draw (-3,\yinit+0.2) -- (10,\yinit+0.2);
		   }
		   \end{scope}
		   
		   \begin{scope}[thick]
			   \draw (-0.4,2.2) -- (0.2,2.2);
			   \draw[<->] (-0.2,0) -- (-0.2,2.2) node [midway,left] {$y_0$};
   		   \draw[blue,very thick] (-3,2.2) -- (10,2.2);
   	   \end{scope}
	   }

	   \draw [->] (-4,0) -- (2.5,0) node [pos=0.9, below] {$x$};
	   \draw [->] (0,-5) -- (0,5) node [pos=0.95, above, sloped] {$y$};

   \end{tikzpicture}
}
\hfill
\subfigure[Zig-zag (case 2).\label{fig:StaticZigZag}]{
   \begin{tikzpicture}[scale=0.75]
	   	   \clip (-1,-1) rectangle ++(6,6);

	   \foreach \oneoverab in {0.75}{
		   \begin{scope}[gray]
		   \foreach \yinit in {-4,-3.5,...,6}{
		      \foreach \cycle in {-4,0,...,8}{
		         \draw (\cycle,\yinit-\oneoverab) cos (\cycle+1,\yinit) sin (\cycle+2,\yinit+\oneoverab) cos (\cycle+3,\yinit) sin (\cycle+4,\yinit-\oneoverab) ;
		      }
		   }
		   \end{scope}
		   
		   \begin{scope}[thick]
		   \foreach \yinit in {2}{
			   \draw (-0.4,\yinit-\oneoverab) -- (0.2,\yinit-\oneoverab);
			   \draw[<->] (-0.2,0) -- (-0.2,\yinit-\oneoverab) node [midway,left] {$y_0$};
			   \draw (-0.4,\yinit+\oneoverab) -- (2.2,\yinit+\oneoverab);
			   \draw[<->] (-0.2,\yinit-\oneoverab) -- (-0.2,\yinit+\oneoverab) node [midway,left] {$\frac{2}{ab}$};
			   \draw (4,\yinit-\oneoverab+0.2) -- (4,-0.4);
			   \draw[<->] (0,-0.2) -- (4,-0.2) node [pos=0.5,below] {$\frac{2\pi}{b}$};
		      \foreach \cycle in {-4,0,...,8}{
		         \draw [blue,very thick] (\cycle,\yinit-\oneoverab) cos (\cycle+1,\yinit) sin (\cycle+2,\yinit+\oneoverab) cos (\cycle+3,\yinit) sin (\cycle+4,\yinit-\oneoverab) ;
		      }
			}
   	   \end{scope}
	   }

	   \draw [->] (-4,0) -- (5,0) node [pos=0.95, below] {$x$};
	   \draw [->] (0,-5) -- (0,5) node [pos=0.95, above, sloped] {$y$};

   \end{tikzpicture}
}
\hfill
\subfigure[Spiral (case 2).\label{fig:StaticSpiral}]{
   \begin{tikzpicture}[scale=0.75]
	   	   \clip (-2.75,-3.25) rectangle ++(6.25,6);

   \begin{scope}[gray]     
   \draw [domain=1000:0,variable=\t,smooth,samples=250]
        plot ({0-\t}: {0.005*\t})
        plot ({120-\t}: {0.005*\t})
        plot ({240-\t}: {0.005*\t});
   \end{scope}
        
   \begin{scope}[thick]     
   
   \fill (330:3.155) circle (2pt); node [below,text width=3.5cm] {initial position};
   \draw [<->] (0,0) -- (330:3.155) node [above,pos=0.6,sloped] {$r_0$};
   \draw  (0,0) ++(0:0.7) arc (0:-30:0.7); \node at (-20:0.7) [right] {$\theta_0$};

   \draw [<->] (0,0) -- (40:2.8) node [above,pos=0.6,sloped] {$r$};
   \draw (0,0) ++(0:0.5) arc (0:40:0.5); \node at (25:0.5) [right] {$\theta$};
        
   \draw [<->] (140:0.48) -- (160:2.2) node [above,midway,sloped] {$2 \pi b$};
   
   \draw (1.4,0) -- (-2.0,0) node [left] {$y_c$};
   \draw (0,0) -- (0,-2.45) node[below]  {$x_c$};

   \draw (330:3.155) -- ++(-4.7,0) node [left] {$y_0$};
   \draw (330:3.155) -- ++(0,-0.88) node[below]  {$x_0$};

   \draw [domain=630:0,variable=\t,smooth,samples=250,blue,very thick]
        plot ({240-\t}: {0.005*\t});        

   \end{scope}

	   \draw [->] (-6,-2.25) -- (3.5,-2.25) node [pos=0.98, below] {$x$};
	   \draw [->] (-1.75,-5) -- (-1.75,2.75) node [pos=0.95, above, sloped] {$y$};

   \end{tikzpicture}
}
\caption{Streamlines of the static route choice for the test problems.}
\end{figure}
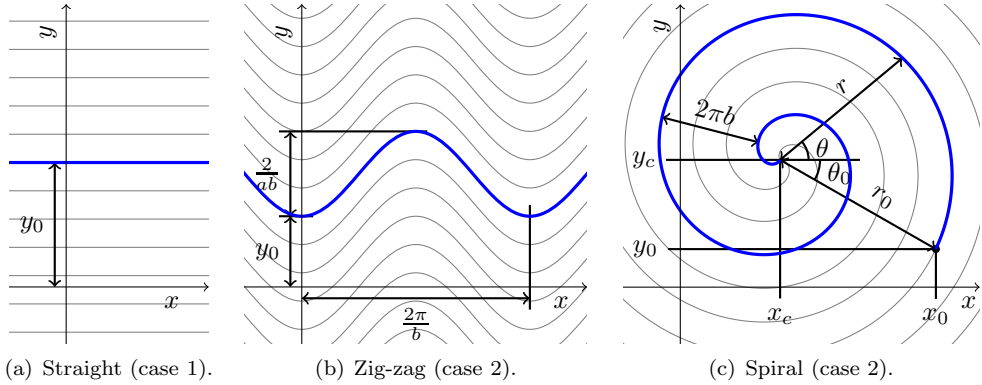

\subsection{Simulation results}
\begin{figure}
\centering
   \subfigure[Densities at times $t=20$s (top) and $t=80$s (bottom).\label{fig:Densities_Straight}]{
   \begin{minipage}[b]{0.46\textwidth}
	\includegraphics[width=1\textwidth]{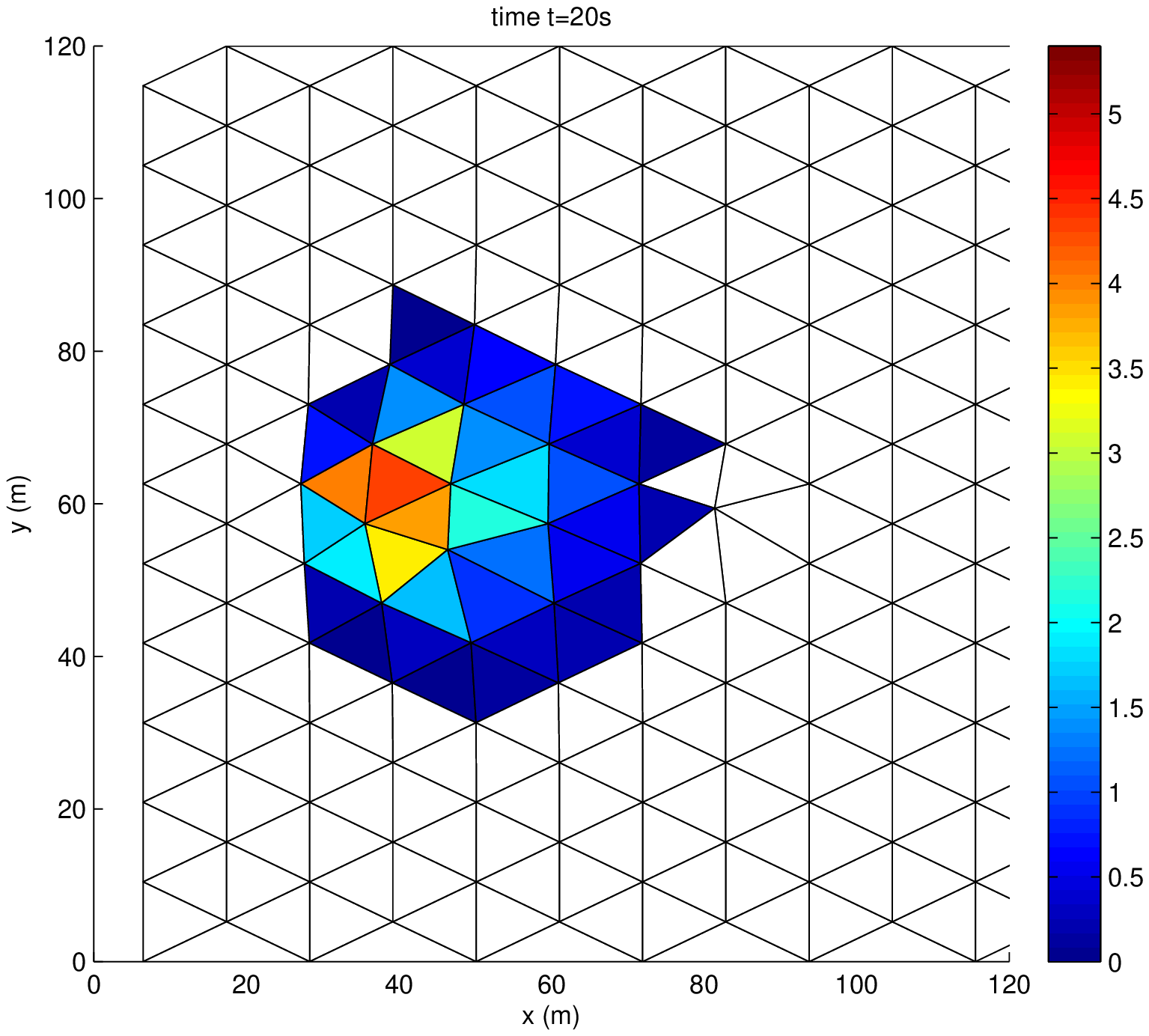}
	\\
	\includegraphics[width=1\textwidth]{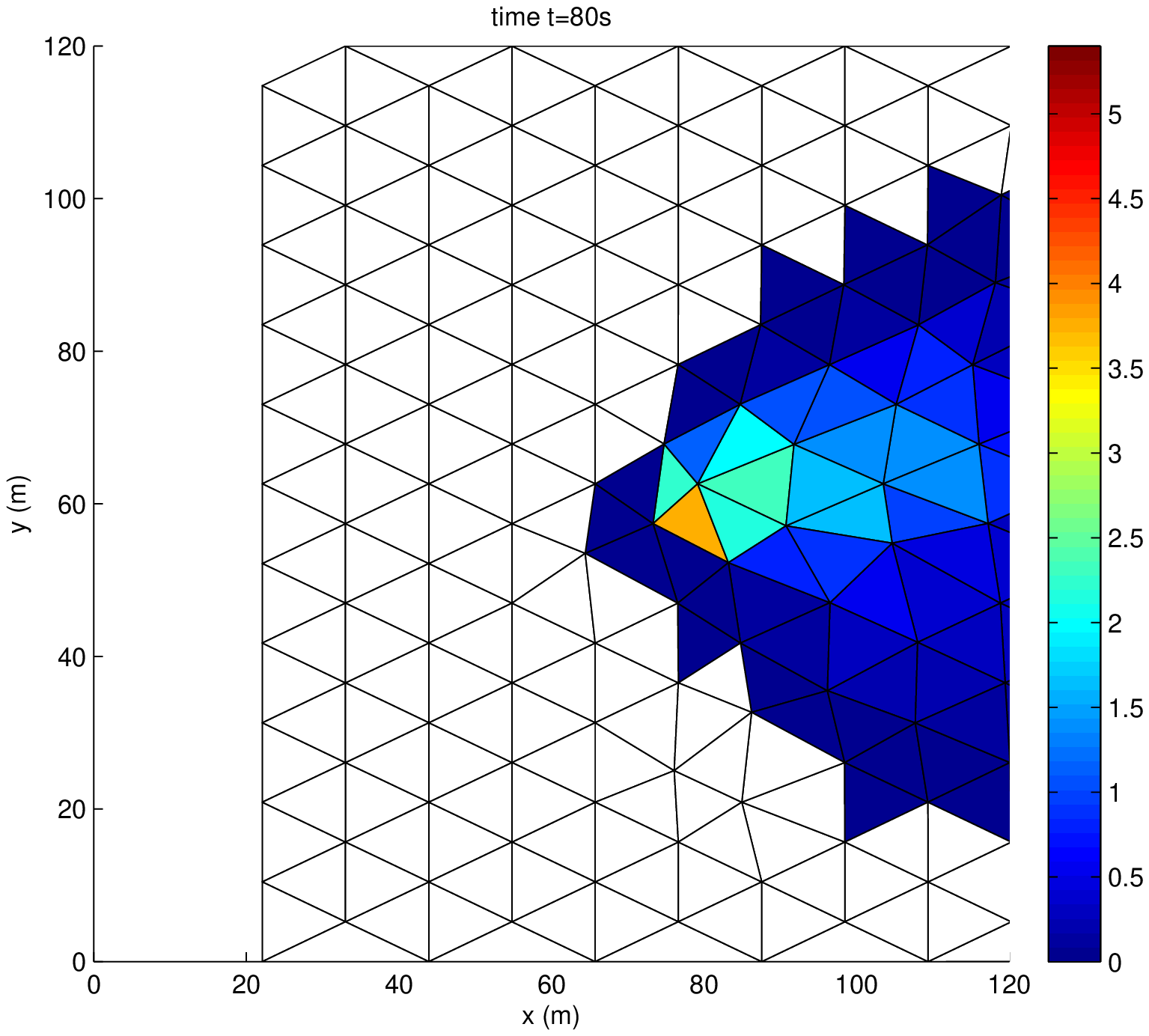}
	\end{minipage}
	}
   \hfill
   \subfigure[Trajectories of centres of grid cells with at least 1 pedestrian in them. Line width corresponds to the number of pedestrians in the grid cell and color indicates the time, both as shown in the legend.\label{fig:Trajectories_legend}]{
   \begin{minipage}[b]{0.46\textwidth}
   \includegraphics[width=1\textwidth,trim=0mm 0mm 0mm 6mm, clip=true]{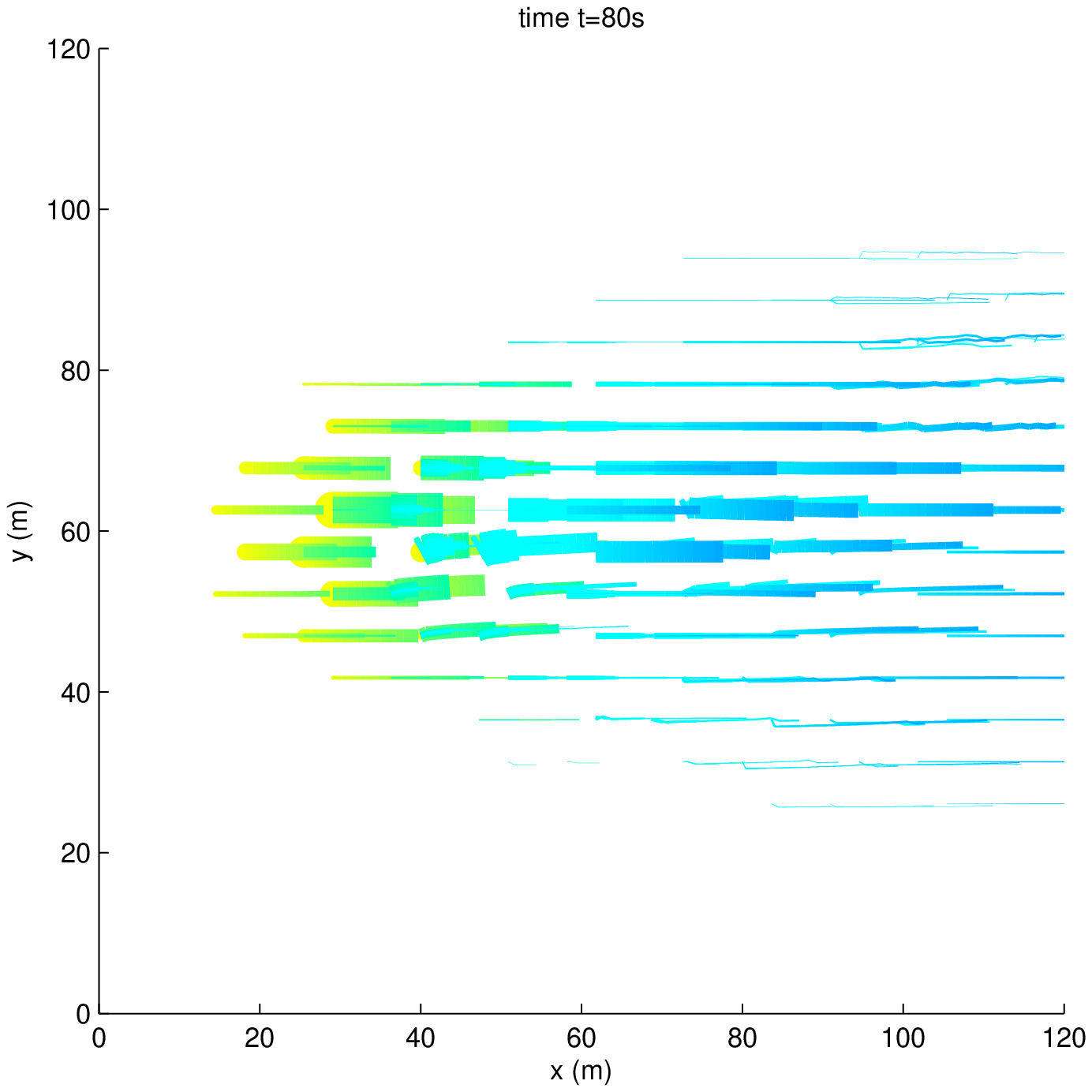}
   \bigskip \\
   number of pedestrians:\\
	\begin{tikzpicture}[scale=0.4]
	\fill (0,0.1pt) -- (0,-0.1pt) -- (12,-2pt) -- (12,2pt) -- cycle;
	\draw (0,4pt) -- (0,-4pt) node [below] {1};
	\draw (4,4pt) -- (4,-4pt) node [below] {50};
	\draw (8,4pt) -- (8,-4pt) node [below] {100};
	\draw (12,4pt) -- (12,-4pt) node [below] {150};
	\end{tikzpicture} 
   \bigskip \\
   time scale (in seconds): \\
	\begin{tikzpicture}[scale=0.4]
	\path[shade, left color=yellow, right color=green] (0,0) rectangle (4,0.3);
	\path[shade, left color=green, right color=cyan]   (4,0) rectangle (8,0.3);
	\path[shade, left color=cyan, right color=blue]    (8,0) rectangle (12,0.3);
	\draw (0,0.3) -- (0,0) node [below] {0};
	\draw (3,0.3) -- (3,0) node [below] {20};
	\draw (6,0.3) -- (6,0) node [below] {40};
	\draw (9,0.3) -- (9,0) node [below] {60};
	\draw (12,0.3) -- (12,0) node [below] {80};
	\end{tikzpicture}
	\end{minipage}
	}
\caption{Densities and trajectories of the centres of the grid cells, for test case 1 (`straight').}
\label{fig:Densities_Trajectories_Straight}
\end{figure}

\begin{figure}
\centering
   \subfigure[Densities at times $t=20$s (top) and $t=80$s (bottom).\label{fig:Densities_ZigZag}]{
   \begin{minipage}[b]{0.46\textwidth}
	\includegraphics[width=1\textwidth]{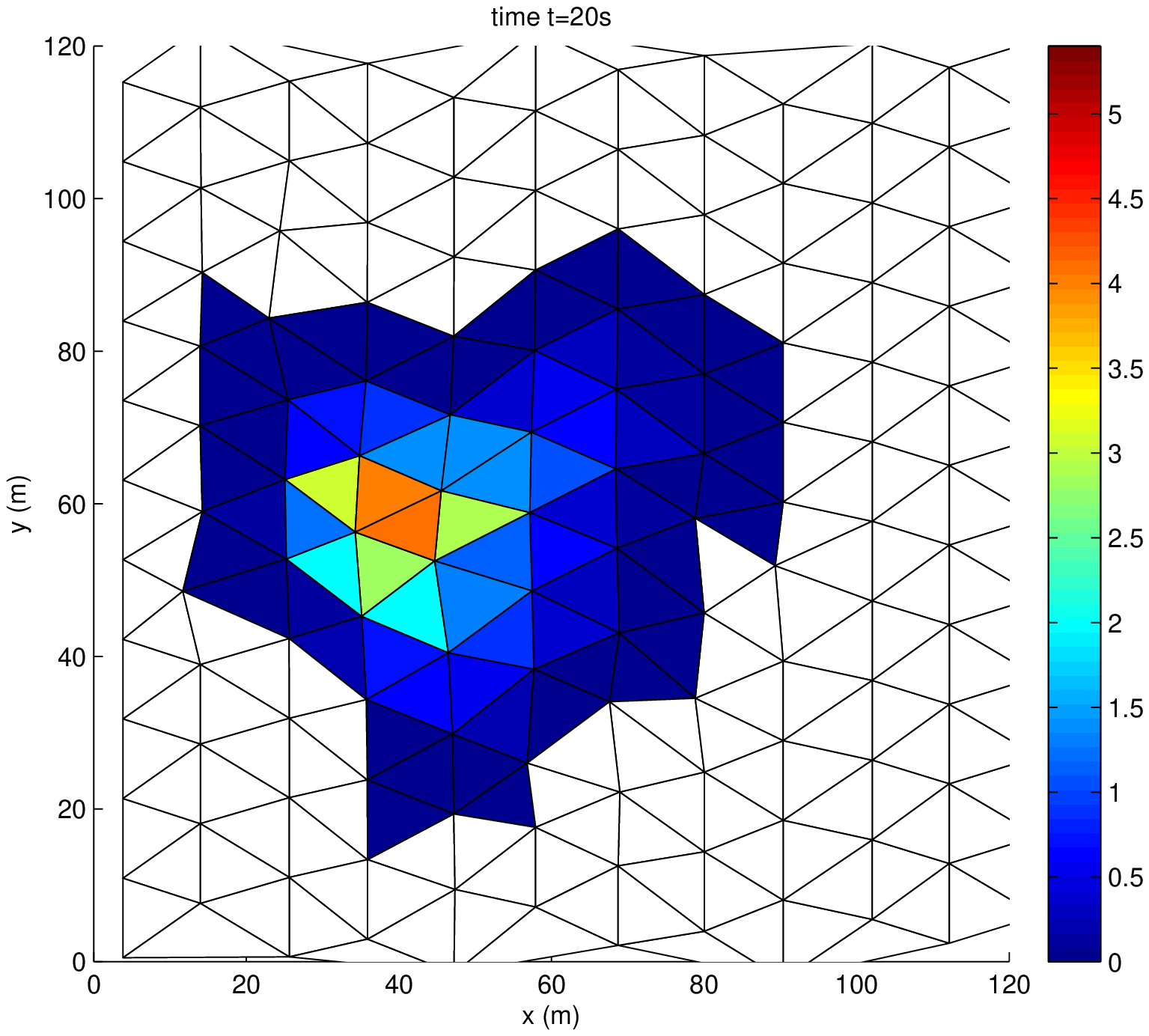}
	\\
	\includegraphics[width=1\textwidth]{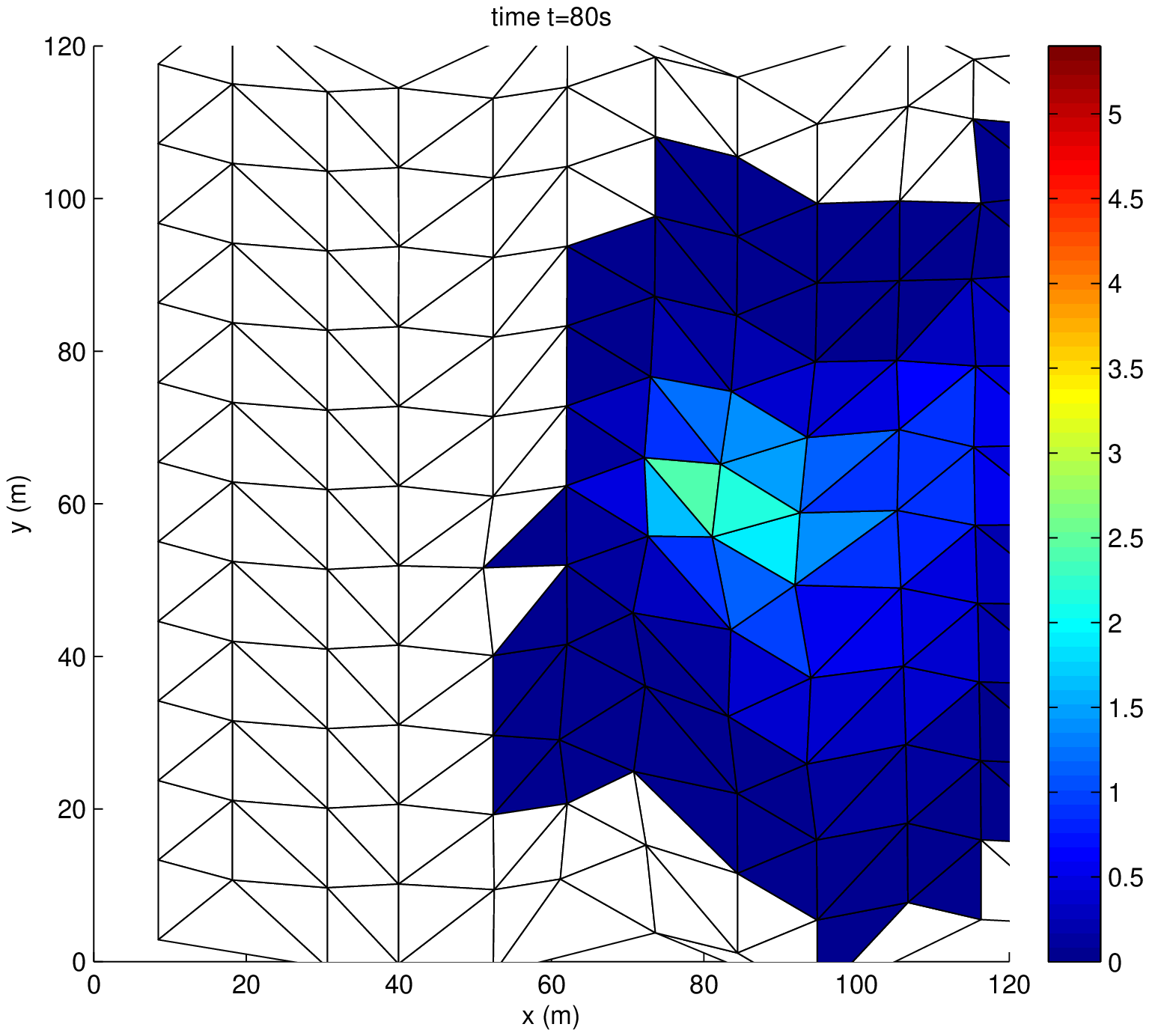}
	\end{minipage}
	}
   \hfill
   \subfigure[Trajectories of centres of grid cells with at least 1 pedestrian in them. Line width corresponds to the number of pedestrians in the grid cell and color indicates the time, both as shown in the legend.\label{fig:Trajectories_legend}]{
   \begin{minipage}[b]{0.46\textwidth}
   \includegraphics[width=1\textwidth,trim=0mm 0mm 0mm 6mm, clip=true]{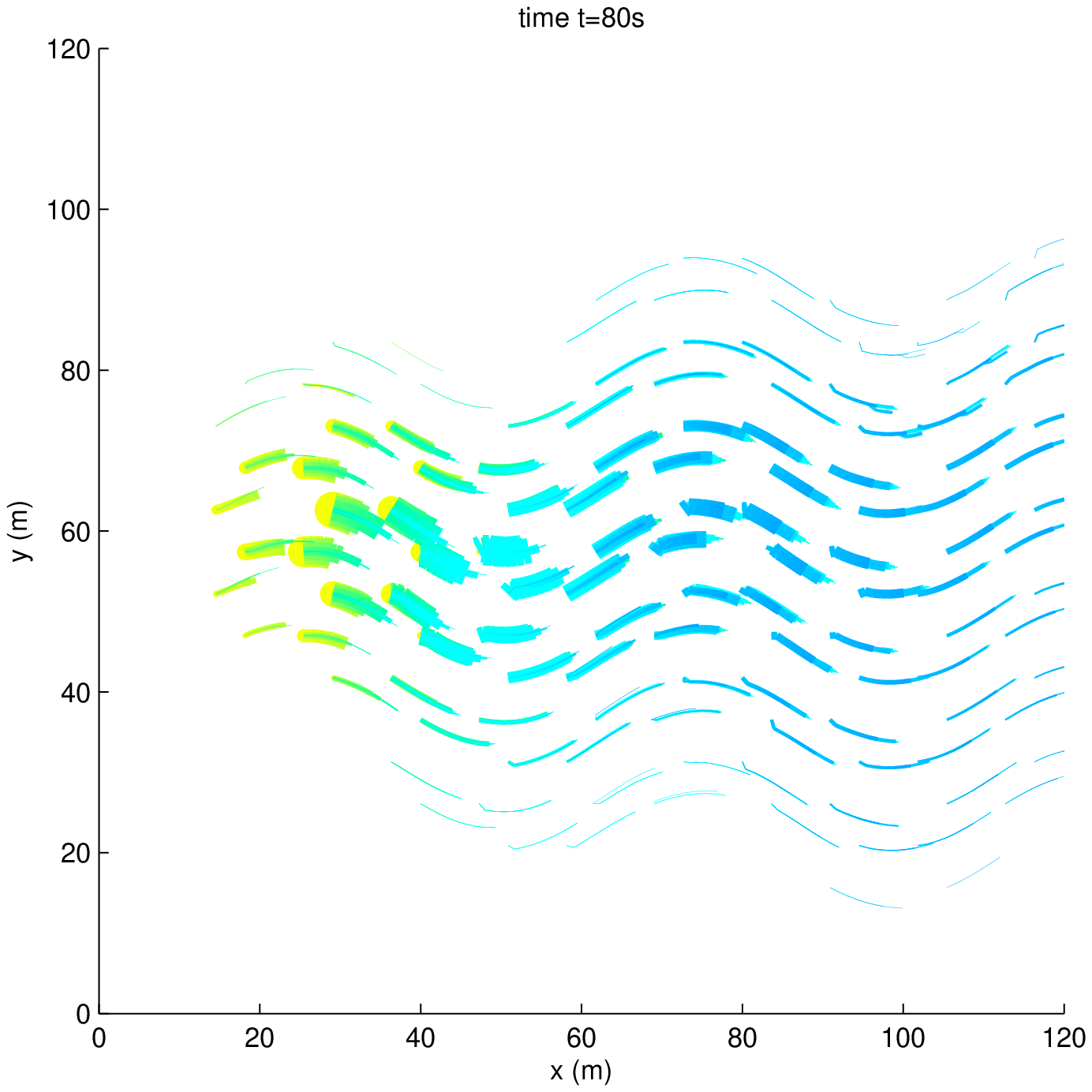}
   \bigskip \\
   number of pedestrians:\\
	\begin{tikzpicture}[scale=0.4]
	\fill (0,0.1pt) -- (0,-0.1pt) -- (12,-2pt) -- (12,2pt) -- cycle;
	\draw (0,4pt) -- (0,-4pt) node [below] {1};
	\draw (4,4pt) -- (4,-4pt) node [below] {50};
	\draw (8,4pt) -- (8,-4pt) node [below] {100};
	\draw (12,4pt) -- (12,-4pt) node [below] {150};
	\end{tikzpicture} 
   \bigskip \\
   time scale (in seconds): \\
	\begin{tikzpicture}[scale=0.4]
	\path[shade, left color=yellow, right color=green] (0,0) rectangle (4,0.3);
	\path[shade, left color=green, right color=cyan]   (4,0) rectangle (8,0.3);
	\path[shade, left color=cyan, right color=blue]    (8,0) rectangle (12,0.3);
	\draw (0,0.3) -- (0,0) node [below] {0};
	\draw (3,0.3) -- (3,0) node [below] {20};
	\draw (6,0.3) -- (6,0) node [below] {40};
	\draw (9,0.3) -- (9,0) node [below] {60};
	\draw (12,0.3) -- (12,0) node [below] {80};
	\end{tikzpicture}
	\end{minipage}
	}
\caption{Densities and trajectories of the centres of the grid cells, for test case 2 (`zig-zag').}
\label{fig:Densities_Trajectories_ZigZag}
\end{figure}

\begin{figure}
\centering
   \subfigure[Densities at times $t=50$s (top) and $t=200$s (bottom).\label{fig:Densities_Spiral}]{
   \begin{minipage}[b]{0.46\textwidth}
	\includegraphics[width=1\textwidth]{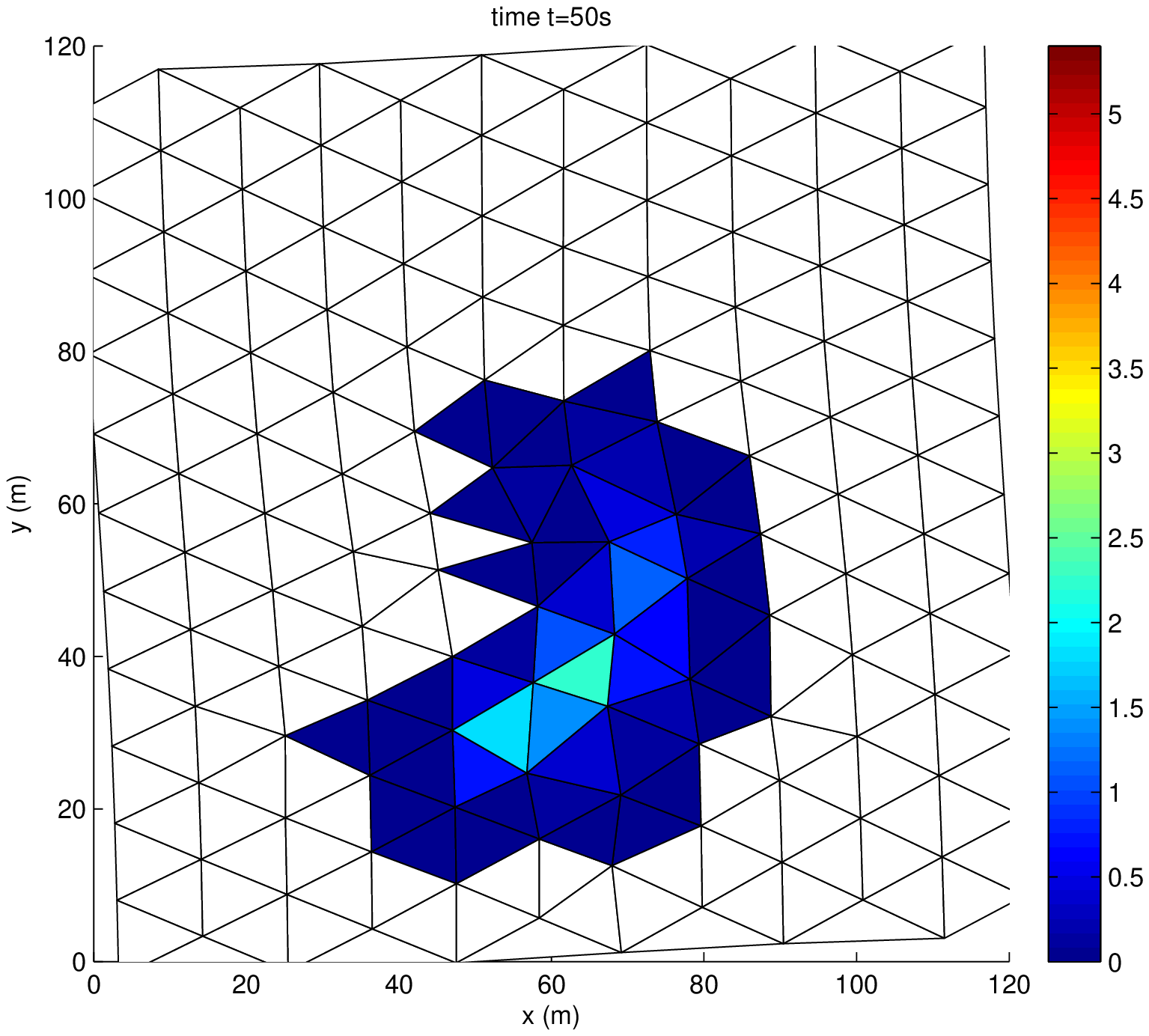}
	\\
	\includegraphics[width=1\textwidth]{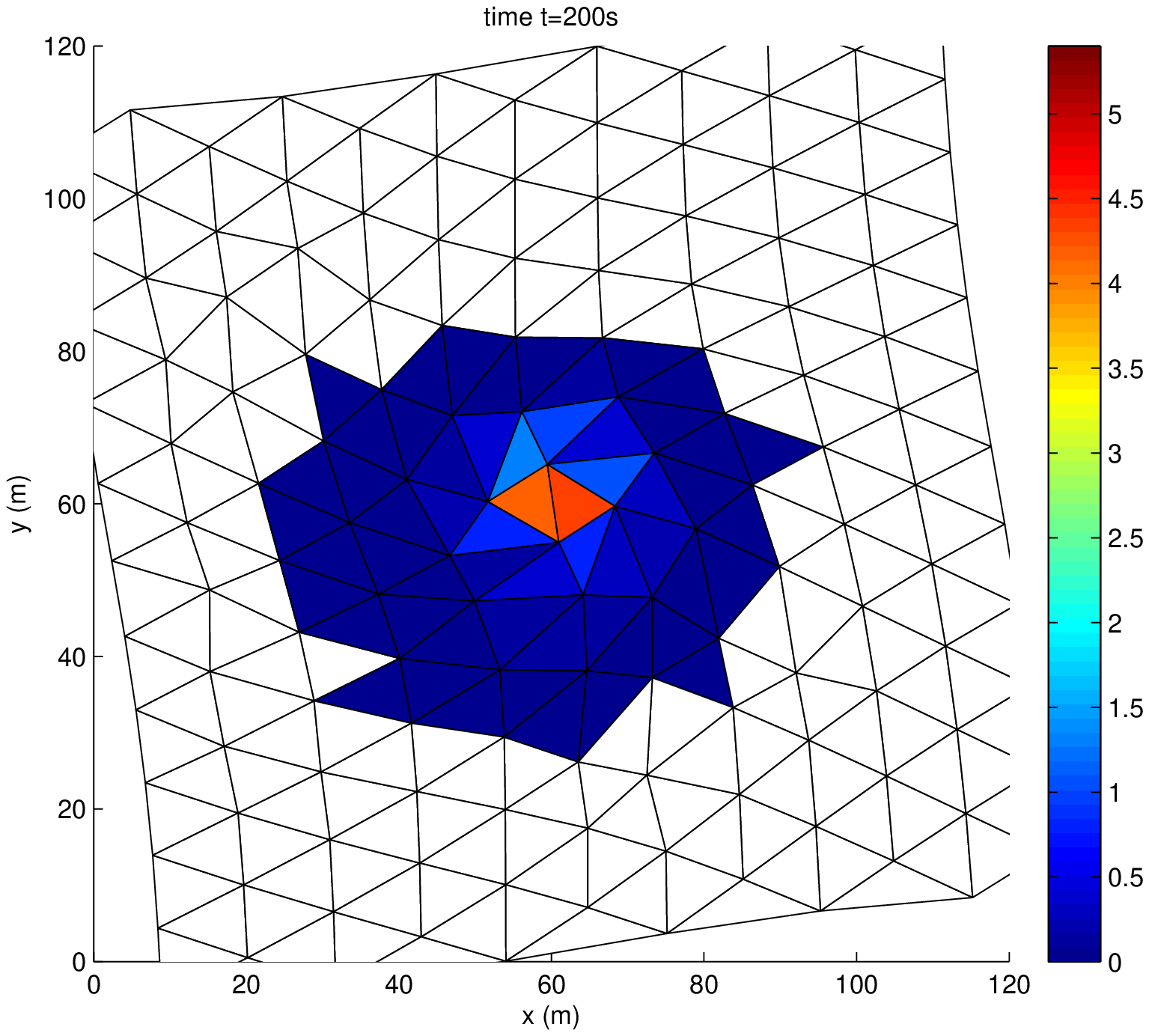}
	\end{minipage}
	}
   \hfill
   \subfigure[Trajectories of centres of grid cells with at least 1 pedestrian in them. Line width corresponds to the number of pedestrians in the grid cell and color indicates the time, both as shown in the legend.\label{fig:Trajectories_legend}]{
   \begin{minipage}[b]{0.46\textwidth}
   \includegraphics[width=1\textwidth,trim=0mm 0mm 0mm 6mm, clip=true]{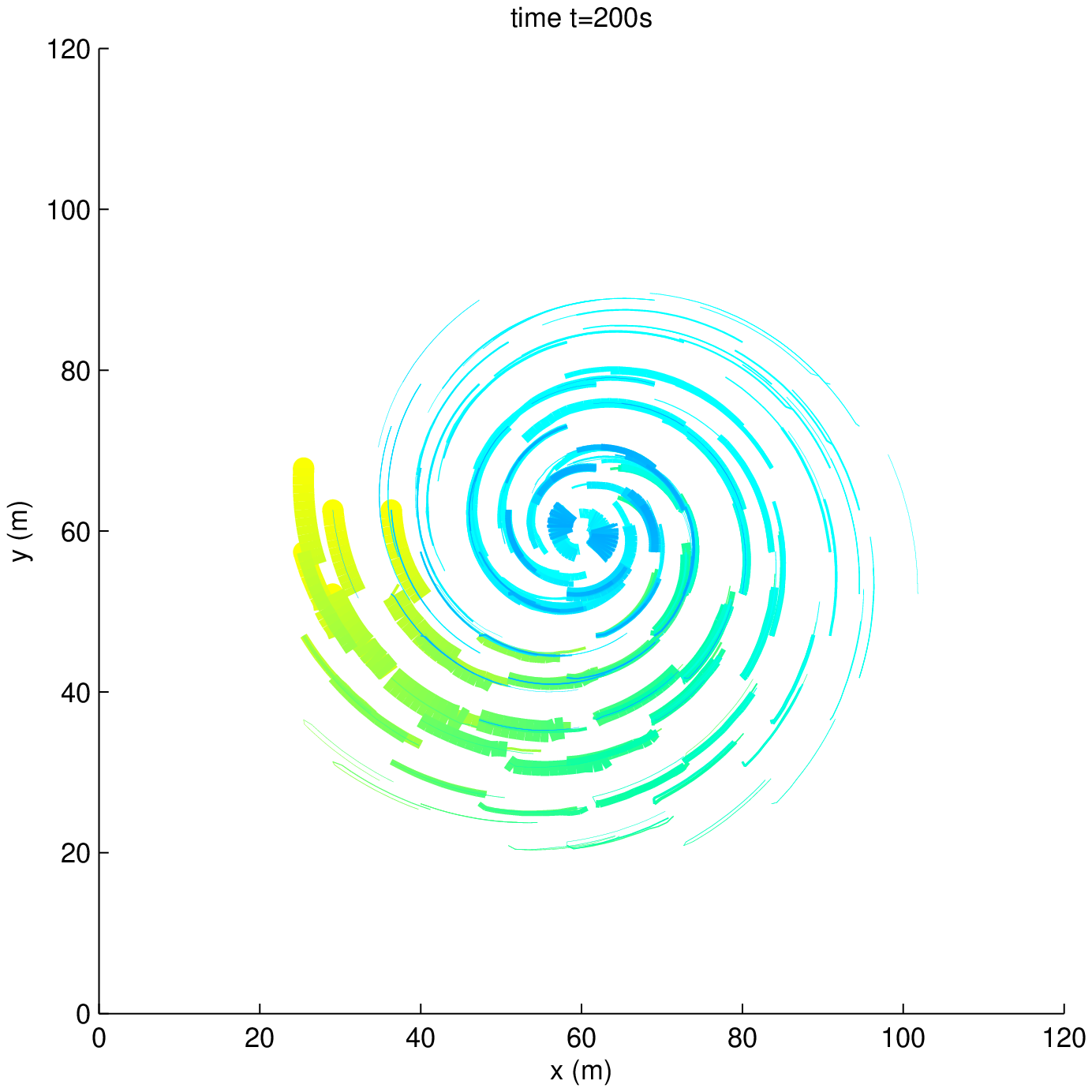}
   \bigskip \\
   number of pedestrians:\\
	\begin{tikzpicture}[scale=0.4]
	\fill (0,0.1pt) -- (0,-0.1pt) -- (12,-2pt) -- (12,2pt) -- cycle;
	\draw (0,4pt) -- (0,-4pt) node [below] {1};
	\draw (4,4pt) -- (4,-4pt) node [below] {50};
	\draw (8,4pt) -- (8,-4pt) node [below] {100};
	\draw (12,4pt) -- (12,-4pt) node [below] {150};
	\end{tikzpicture} 
   \bigskip \\
   time scale (in seconds) \\
	\begin{tikzpicture}[scale=0.4]
	\path[shade, left color=yellow, right color=green] (0,0) rectangle (4,0.3);
	\path[shade, left color=green, right color=cyan]   (4,0) rectangle (8,0.3);
	\path[shade, left color=cyan, right color=blue]    (8,0) rectangle (12,0.3);
	\draw (0,0.3) -- (0,0) node [below] {0};
	\draw (3,0.3) -- (3,0) node [below] {50};
	\draw (6,0.3) -- (6,0) node [below] {100};
	\draw (9,0.3) -- (9,0) node [below] {150};
	\draw (12,0.3) -- (12,0) node [below] {200};
	\end{tikzpicture}
	\end{minipage}
	}
\caption{Densities and trajectories of the centres of the grid cells, for test case 3 (`spiral').}
\label{fig:Densities_Trajectories_Spiral}
\end{figure}

Figures \ref{fig:Densities_Trajectories_Straight}, \ref{fig:Densities_Trajectories_ZigZag} and \ref{fig:Densities_Trajectories_Spiral} show the simulation test results. It shows that the pedestrians are directed by the static route choice, while also being repelled by other pedestrians. The results look plausible, but will be analysed more in-depth in future studies.

%

\section{Conclusion and future work}
\label{sec:Conclusion}
We have introduced a numerical method for crowd flow simulation in two dimensions, based on the Lagrangian coordinate system. We have shown that the results are plausible. In the development of the method, we made some choices which could be debated. In future research, we plan to look into the choices and alternatives in more detail. For example, we consider using grid cells shaped as general polygons, not necessarily triangles and we will consider alternative remeshing procedures, including local remeshing. We will study both the accuracy and the computation time of the newly developed methods in more depth. Furthermore, we will compare the results with those of alternative computational methods, such as those based on the Eulerian coordinate system.

Future extensions of the method include the introduction of boundary conditions and multiple groups that have different walking or route choice characteristics.

\bibliographystyle{plainnat}
\bibliography{Lagrange2D}

\end{document}